\documentclass{article}
\usepackage{amsfonts}
\usepackage{amsmath}
\usepackage[T1]{fontenc}
\usepackage[utf8]{inputenc}
\usepackage{authblk}
\usepackage{setspace}
\usepackage{cite}
\usepackage{amsthm}
\usepackage{CJK,fancybox,amssymb,enumerate}
\usepackage{color}

\usepackage{graphicx}
\usepackage{float}
\usepackage{subfig}
\usepackage{indentfirst}
\usepackage{epstopdf}
\usepackage{xr}
\usepackage{hyperref}
\usepackage{cleveref}
\numberwithin{equation}{section}
\newtheorem{definition}{Definition} [section]
\newtheorem{theorem}{Theorem}[section]
\newtheorem{lemma}{Lemma}[section]

\newtheorem{corollary}{Corollary}[section]

\date{}

\title{Dynamical behavior of a colony migration system: Do colony size and quorum threshold affect collective-decision?
\footnote{This work was partially supported by NSF-DMS (Award Number 1716802\&2052820);  NSF- IOS/DMS (Award Number 1558127), DARPA-SBIR 2016.2 SB162-005, and The James S. McDonnell Foundation 21st Century Science Initiative in Studying Complex Systems Scholar Award (UHC Scholar Award 220020472). The work of the first author was partially supported by the Scholarship Foundation of China Scholarship Council (award 201906840071). The work of the second author was supported by the National Natural Science Foundation of China through grants 12071217 and 11671206.}}

\author[1]{Lisha Wang\thanks{E-mail address:\ lswang.math@gmail.com
}}
\author[2]{Zhipeng Qiu\thanks{E-mail address:\ nustqzp@njust.edu.cn
}}
\author[3]{Takao Sasaki\thanks{E-mail address:\ takao.sasaki@uga.edu
}}
\author[4]{Yun Kang\thanks{Corresponding author. E-mail address:\ yun.kang@asu.edu}}

\affil[1] {\small Department of Mathematics, Nanjing University of Science and Technology, Nanjing 210094,
	China}

\affil[2]{\small Center for Basic Teaching and Experiment, Nanjing University of Science and Technology, Jiangyin 214443,
	China}

\affil[3]{\small Odum School of Ecology, University of Georgia, Athens, GA 30602 USA}

\affil[4]{\small Sciences and Mathematics Faculty, College of Integrative Sciences and Arts, Arizona State University, Mesa, AZ 85212 USA}

\begin{document}
\maketitle

{\bf Abstract.}Social insects are ecologically and evolutionarily most successful organisms on earth, which can achieve robust collective behaviors through local interactions among group members. Colony migration has been considered as a leading example of collective decision-making in social insects. In this paper, a piecewise colony migration system with recruitment switching is proposed to explore underlying mechanisms and synergistic effects of colony size and quorum on the outcomes of collective decision. The completed dynamical analysis for the non-smooth system (including the dynamics on subsystems, switching surface, and full system) is performed, and the sufficient conditions significantly related to colony size for the stability of equilibria are provided. The theoretical results suggest that large colonies are more likely to emigrate to a new site. More interesting findings include but not limit to: (a) the system may exhibit oscillation when the colony size is below a critical level; (b) the system may also exhibit a bistable state, and colonies migrate to a new site or the old nest depending on their initial sizes of recruiters. Bifurcation analysis shows that the variations of colony size and quorum threshold greatly impact the dynamics. The results suggest that it is important to distinguish two populations of recruiters in modeling. This work may provide important insights on how simple and local interactions achieve the collective migrating activity in social insects. 


\section{Introduction}
Social insects have been studied extensively since they are typical group living organisms with collective decision-making behaviors \cite{camazine2020self,sumpter2010collective,seeley2010honeybee,sasaki2012groups,seeley2009wisdom,feng2021recruitment}. The members in these groups can make a colony-level choice by individual communication and acting with decision rules to seek a consensus outcome \cite{pratt2002quorum}. Without any central control, their contribution to the global decision only comes from the local and limited information resource \cite{pratt2005agent}. However, it still enables the colony an accurate and efficient decision from complex environment. Social insects with collective decision-making behavior range from the foraging honey-bee to the migrating ants, all of which can perform complex organizational activities without well-informed leaders\cite{seeley1999group,visscher1999collective,mallon2001individual,bonabeau1997self,pratt1998decentralized,hirsh2001distributed}. These biological {phenomena} encourage more perspectives on understanding the relationship between individual behavioral rules and the overall ability behind performing complex activities.

Colony migration of social insects is {one of} leading example of collective decision-making behavior \cite{pratt2005behavioral}. The  colony as a whole can move to a suitable nest rather than splitting. {They typically} achieve this {consensus decision through} a high degree of communication and coordination among group members \cite{pratt2002quorum}. Colony migration in the genus \textit{Temnothorax} (formerly \textit{Leptothorax}) is a particularly promising subject. \textit{Temnothorax} {ants typically} live in rock crevices and are likely to require migrating frequently {due to fragility of their nest sites} \cite{partridge1997field}. In the laboratory,  ants can be easily marked and monitored by taking advantage of their small colony size (usually a few hundred workers). {Using these detailed individual data, extensive investigations \cite{pratt2005agent,beckers1989colony,moglich1978social,beckers1990collective} have revealed underlying processes during a migration} in this genus. Generally, migrations are {initiated} only by active workers, about one-third of the colony, who search for potential new homes, assess their quality, and recruit nestmates to the finds. Understanding this emergence of colony migration would provide insights into the study of a wide array of collective decision-making behavior.

Recently, increasing experimental work has promoted a deeper understanding on the process of colony migration incorporating complicated individual behavior and decision rules \cite{todd2000precis,dornhaus2004ants}. Mallon et al. \cite{mallon2001individual} showed that the ants may contribute to the collective decision through quality-dependent difference of recruitment latency, i.e., the individuals take less time to initiate recruitment to a superior than to a mediocre site. Pratt et al. \cite{pratt2002quorum} found that the scout assesses new finding site and then recruit nestmates through tandem running until a quorum threshold is reached, at which point the ant changes to transporting nestmates. In \cite{pratt2005quorum}, Pratt revealed that the ants measure the achievement of a quorum through their rate of direct encounters with nestmates. Sasaki et al. \cite{sasaki2019rational} studied the rationality of time investment during nest-site choice, and the results show that the isolated ants took more time to complete the migration when choosing between two similar nests, but the whole colonies rationally made faster decisions. {These experimental works exhibit extensive interesting phenomena of emergence of collective decision making from individuals.} Thus, to full understand the collective decision-making in colony migration, it is necessary to investigate the mechanism underlying it.

Mathematical model is the powerful tool to gain insights into deeper analysis on the colony mechanism and to explain collective performance in migration. Most of recent mathematical works of colony migration concentrate on simulating colony-wide trends by using agent based model \cite{pratt2005agent,pratt2006tunable,tofts1992describing,de1998modelling,sumpter2001ants}. However, it is also necessary to develop models for analyzing the dynamics and generating testable predictions in changing environment. Differential equations can provide a better understanding of how multiple components of colony migration interact with each other. In \cite{pratt2002quorum}, Pratt et al.  have proposed differential equations to explore how a quorum can help colonies choose between two sites with different quality, and the simulations show that the colony splits into different sites when the quorum is too small and reach a consensus on nest-choice by increasing quorum threshold. Assis et al. \cite{assis2009decision} have presented a differential model to describe the competition of the different sites, and they clarify that the threshold factor and the flux of resource provided by the colony play roles in decision-making. {Although some agent based models and differential equations models have been proposed to explore the colony migration behavior in social insects, it is still in an early stage to rigorously analyze the collective migration process by using mathematical tools. Motivated by  \cite{pratt2002quorum} and the recent work in \cite{pratt2005agent}, we develop an ODE model that incorporates complicated migration rules and provide some biological implications from novel interesting mathematical studies.}

Increased evidence suggested that the variation of colony size significantly affects collective behaviors in social insects. Many works have shown a positive correlation between group size and information flow rate \cite{burkhardt1998individual,karsai1998productivity,gordon1999encounter}. Larger colony size may display a higher level of division of labor and allocation of tasks \cite{gautrais2002emergent,feng2021dynamics}, more effective exploration with lower risk aversion \cite{dornhaus2006colony,herbers1981reliability},  and can better resist random disturbance of local information acquisition \cite{mallon2001individual}. In some cases, the colony size can also affect the time needed to make a decision and the methods used in recruitment in group activities \cite{beckers1989colony,planque2010recruitment}.
Dornhaus et al. \cite{dornhaus2006colony} studied the influence of colony size on collective decision-making in the colony migration. The results show that the quorum threshold may remain constant with the size of natural colonies or be proportional to the size of manipulated colonies. All the biological observation support the hypothesis that colony size is important to collective decision-making in ants. Hence, it is also necessary to evaluate the potential impact of colony size as well as the synergistic effect of colony size and quorum threshold on the outcomes of migrations. In this paper, we develop a mathematical model to describe the process of colony migration in dynamical environment. Our proposed model is expected to address the following ecological questions in social insects from our mathematical studies:
\begin{itemize}
	\item How does the colony size affect the migration result?
	\item  How do synergies of colony size and quorum threshold regulate migration dynamic behaviors?
\end{itemize}

The structure of this article is organized as follows: In Section \ref{Section2}, we provide the biological background of colony migration and derive a migrating system described by piecewise differential equations. In Section \ref{Section3}, we perform the mathematical analysis of our model. In Section \ref{Section4}, we classify the dynamical behaviors of the colony migration system. In Section \ref{Section5}, we investigate the synergistic effects of colony size and quorum threshold on the dynamics of system through bifurcation analysis. In Section \ref{Section6}, we provide a conclusion of our results and the potential outlook of our current work. 

\section{Model derivations}
\label{Section2}


We start with a simple description of workers's behavior during migration processes. Generally, all active workers follow a strategy of graded commitment to the site they have found, with transitions to higher levels depending both on the quality of site and on the interactions among nestmates \cite{pratt2005agent}. At the lowest level of commitment, the searchers enter new finding site and stay inside for an independent assessment. The duration of assessment is inversely related to the quality of new site. At the next level, the workers start to recruit other active workers via tandem runs, in which a single follower is led from the old nest to the new site. The new arrivals would make their independent decisions about whether to recruit. Once the number of active workers presented at new site reaches a quorum threshold, the workers enter the highest level of commitment. They carry remaining nestmates and brood items to new home by transportation recruitment. At any level of commitment,  workers may leave the current site with a probability and search the surrounding area again for a new potential site.
\begin{figure}[htbp]
	\centering
	\includegraphics[width=8.5cm]{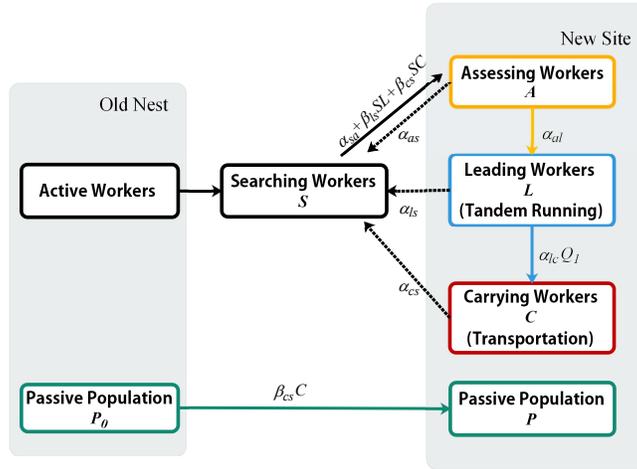}
	\caption{Model diagram of single-nest colony migration.}
	\label{ModelDiagram1}
\end{figure}

The model presented in this paper is based on assumed processes showed in Figure \ref{ModelDiagram1}. We consider the most typical scenario of emigration, namely only one potential site is available near the old nest. Assume that the colony in the old nest has a total of $N$ workers where  $\rho N$ is  active population, and $(1-\rho)N$ is passive population. According to the biological description, each active worker should be in one of the following four classes: searching workers denoted by $S$, assessing workers denoted by $A$, leading workers denoted by $L$, and carrying workers denoted by $C$. The passive population remaining in the old nest is denoted as $P_0$, and the passive population moving to the new site is denoted as $P$. A transition diagram between different classes of populations is depicted in Figure \ref{ModelDiagram1} whose assumptions are showed below:
\begin{enumerate}[(a)]
	\item  During colony migration, the total number of workers in this colony is constant, i.e., $N=P_0+P+S+A+L+C$.
	
	\item \textbf{Searching workers $S$.} The number of searching workers $S$ depends on the rates at which assessing workers $A$, leading workers $L$ and carrying workers $C$ join searching workers $S$, $\alpha_{as}A$, $\alpha_{ls}L$ and $\alpha_{cs}C$ respectively; the rate at which searching workers independently find new site and join assessing workers $A$, $\alpha_{sa}S$; the rate at which the searching workers transit into assessing workers $A$ by interaction with leaders, $\beta_{ls}SL$; the rate at which searching workers transit into assessing workers $A$ by interaction with carriers, $\beta_{cs}SC$. Therefore, the population dynamics of the searching workers $S$ could be described by: 
	\begin{displaymath}
	\begin{aligned}
	\frac{dS}{dt}=&-\underbrace{\beta_{cs} SC}_{\text{$S$ transits to $A$ after interaction with $L$}}-\underbrace{\beta_{ls}SL}_{\text{$S$ transits to $A$ after interaction with $C$}}\\
	&-\underbrace{\alpha_{sa}S}_{\text{$S$ transits to $A$ independently}}+\underbrace{\alpha _{as}A+\alpha _{ls}L+\alpha _{cs}C}_{\text{ The transition from $A$,$L$ and $C$ to $S$}}.
	\end{aligned}
	\end{displaymath}
	\item \textbf{Assessing workers $A$.} The number of assessing workers $A$ depends on the rate at which searching workers $S$ join assessing workers $A$ through independently finding a new site, $\alpha_{sa}S$;  the rates at which the searching workers $S$ transit into assessing workers $A$ by interactions with leaders and carriers respectively, $\beta_{ls}SL$ and $\beta_{cs}SC$; the rate at which assessing workers $A$ join searching workers $S$, $\alpha_{as}A$; the rate at which assessing workers $A$ join leading workers $L$, $\alpha_{al}A$.	Therefore, the population dynamics of the assessing workers $A$ could be described by:
	\begin{displaymath}
	\begin{aligned}
	\frac{dA}{dt}=&\underbrace{\beta_{ls}SL}_{\text{$S$ transits to $A$ after interaction with $L$}}+\underbrace{\beta_{cs} SC}_{\text{$S$ transits to $A$ after interaction with $C$}}\\
	& +\underbrace{\alpha_{sa} S}_{\text{$S$ transits to $A$ independently}}-\underbrace{\alpha _{as}A}_{\text{The transition to $S$}}-\underbrace{\alpha_{al}  A}_{\text{The transition to $L$}}.
	\end{aligned}
	\end{displaymath}
	
	\item \textbf{Leading workers $L$.} The number of leading workers $L$ depends on the rate at which assessing workers $A$ join leading workers $L$, $\alpha_{al}A$; the rate at which leading workers $L$ join to searching workers $S$, $\alpha_{ls}L$; the rate at which leading workers $L$ join carrying workers $C$, $\alpha_{lc}Q_1L$ where $Q_1$ is the probability of switching recruitment decision. The recruitment decision is scored as either $0$ or $1$ depending on the size between total active workers at new site  ($A+L+C$) and quorum threshold ($\Theta$). Specifically, $Q_1=1$ if $A+L+C> \Theta$, and $Q_1=0$ if $A+L+C< \Theta$. Therefore, the population dynamics of the leading workers $L$ could be described by:
	\begin{displaymath}
	\frac{dL}{dt}=\underbrace{\alpha_{al}  A}_{\text{{The transition from $A$}}}-\underbrace{\alpha_{lc} Q_1 L}_{\text{{Recruitment switching}}}-\underbrace{\alpha _{ls}L}_{\text{$L$ transits to $S$}}.
	\end{displaymath}
	
	\item \textbf{Carrying workers $C$.} The number of carrying workers $C$ depends on the rate at which leading workers $L$ join carrying workers $C$, $\alpha_{lc}Q_1L$; the rate at which carrying workers $C$ join searching workers $S$, $\alpha_{cs}C$.	Therefore, the population dynamics of the carrying workers $C$ could be described by:
	\begin{displaymath}
	\frac{dC}{dt}=\underbrace{\alpha_{lc} Q_1 L}_{\text{{Recruitment switching}}}-\underbrace{\alpha _{cs}C}_{\text{$C$ transits to $S$}}.
	\end{displaymath}
	\item \textbf{Passive population $P$ at the new site.} The size of passive population $P$ depends on the rate at which the passive workers $P$ is transported from old nest to new site by carrying workers $C$, $\beta_{cs} C\left[(1-\rho)N-P\right]$. For single-nest emigration, there is no output of passive population $P$. Therefore, the population dynamics of the passive population $P$ could be described by:
	\begin{displaymath}
	\frac{dP}{dt}=\underbrace{\beta_{cs} C\left[(1-\rho)N-P\right]}_{\text{Passive population is carried from old nest to new site}}.
	\end{displaymath}
\end{enumerate}

Based on the above assumptions, we have the following differential equations to describe the dynamics of colony migration:
\begin{equation}\label{colonymigration}
\begin{aligned}
\frac{dS}{dt}&=-\alpha_{sa} S-\beta_{ls}SL-\beta_{cs} SC+\alpha _{as}A+\alpha_{ls}L+\alpha _{cs}C,\\
\frac{dA}{dt}&=\alpha_{sa} S+\beta_{ls}SL+\beta_{cs} S C-\alpha _{as}A-\alpha_{al}  A,\\
\frac{dL}{dt}&=\alpha_{al}  A-\alpha_{lc} Q_1 L-\alpha _{ls}L,\\
\frac{dC}{dt}&=\alpha_{lc} Q_1 L-\alpha _{cs}C,\\
\frac{dP}{dt}&=\beta_{cs} C\left[(1-\rho)N-P\right],
\end{aligned}
\end{equation}
where $Q_1$ is a switching function defined as follows
\begin{displaymath}
\left\{
\begin{array}{lr}
Q_1=0, \quad \text{if} \quad A+L+C<\Theta, \\
Q_1=1, \quad \text{if} \quad A+L+C>\Theta.\\
\end{array}
\right.
\end{displaymath}
For Model \eqref{colonymigration}, all variables and parameters are listed in Table \ref{table01}. Among these parameters, $\beta_{ls}$ is the recruitment rate by leaders and $\beta_{cs}$ is the recruitment rate by carriers.  For \textit{Temnothorax} ants, recruitment rate of carrier is more rapidly than that of leader, i.e., $\beta_{ls} <\beta_{cs}$. However, there exists opposite situation in other species of social insects, such as \textit{Diacamma indicum} \cite{kaur2017characterization}. Therefore, within the framework of our model in this paper, we also consider the case that $\beta_{ls} \geq \beta_{cs}$.

\textbf{Notes.}  Our work is motivated by the differential equations model in \cite{pratt2002quorum} and the agent based model in \cite{pratt2005agent}. Compared with the model in \cite{pratt2002quorum}, Model \eqref{colonymigration} has three innovations: (i) The model in \cite{pratt2002quorum} incorporates only three types of active populations including searchers, assessors and recruiters, while our model has one more component, i.e., the recruiters are divided into population $L$ and population $C$. (ii) We add the transitions of workers from assessing, leading or carrying population to searching population.  (iii) We assume the nonlinear interactions between pope sure that the available workers $S$ can transit ulation $L$ or $C$ and population $S$ which makinto group $A$ through the physical/signal contacts with leaders or carriers. All these hypotheses are biological relevant. Although the agent based model in \cite{pratt2005agent} includes four component and considers the transitions to searching population, it has not been mathematically analyzed in detail. We incorporate these assumptions into our model to investigate collective migration in social insects by using rigorous mathematical proofs and carefully performed bifurcation analysis.

\begin{table}[tbhp]
	{\footnotesize
		\caption{Descriptions of parameters involved in Model \eqref{colonymigration} and their values taken within the range described in previous literature sources \cite{pratt2005behavioral,pratt2005quorum,pratt2002quorum,pratt2005agent,stroeymeyt2010improving,stroeymeyt2011experience,o2016migration}. }\label{table01}
		\begin{center}
			\begin{tabular}{|c|c|c|c|} \hline%
				\bf Parameter & \bf Description & \bf Units  & Values \\ \hline
				$S$	&  Density of searching population & nbr & -\\
				$A$	& Density of assessor & nbr  & -\\
				$L$	& Density of leader  & nbr  & -\\
				$C$	& Density of carrier  & nbr  & -\\
				$P$	& Density of passive workers at new site  & nbr  & -\\
				$P_0$	& Density of passive population at old nest  & nbr  & -\\
				$N$  &  Total number of workers in colony & nbr  & [0, 350]\\
				$\rho$	&  Proportion of active workers & -  & 0.25\\
				$\alpha_{sa}$	& The discovery rate of new site & $\text{min}^{-1}$  & [0.01, 0.15]\\
				$\alpha_{al} $ & The transition rate from assessors to leaders & $\text{min}^{-1}$  & [0.007, 0.2] \\
				$\alpha_{lc} $ & The transition rate from leaders  to carriers & $\text{min}^{-1}$  & [0.15, 0.28]\\
				$\beta_{ls}$	& The rate at which leaders recruit searchers &  $\text{(min ant)}^{-1}$  & [0.004, 0.049]\\
				$\beta_{cs}$	& The rate at which carries recruit searchers & $\text{(min ant)}^{-1}$  & [0.0025, 0.079]\\
				$\Theta$ &  Quorum threshold  & nbr  & [0, 50]\\
				$\alpha _{as}$  &  The transition rate from assessors to searchers & $\text{min}^{-1}$  & [0.24, 0.5]\\
				$\alpha _{ls}$ & The transition rate from leaders to searchers & $\text{min}^{-1}$  & [0.018, 0.12]\\
				$\alpha _{cs}$ & The transition rate from carriers to searchers & $\text{min}^{-1}$  & [0.05, 0.07]\\
				\hline
			\end{tabular}
		\end{center}
	}
\end{table}


\section{Mathematical Analysis}\label{Section3}
In this section, we perform mathematical analysis on the existence and stability of equilibria of the colony migration model \eqref{colonymigration}.  Let $\beta=\max\{\beta_{ls}, \beta_{cs}\}$ be the maximum recruitment rate of the colony and $\sigma=\min\{\alpha_{as}, \alpha_{ls}, \alpha_{cs}\}$ be the minimum transition rate from other groups to the searching group $S$. The basic dynamical result regarding Model \eqref{colonymigration} is shown below.

\begin{theorem}\label{positiveinvariant}
	Model \eqref{colonymigration} is positive invariant in $\mathbb{R}_{+}^{5}$, and every trajectory of Model \eqref{colonymigration} attracts to the compact set 
	\begin{displaymath}
	\Omega=\left\{(S, A, L, C, P)\in\mathbb{R}_{+}^{5}:S+A+L+C=\rho N, 0\leq P\leq (1-\rho)N \right\}
	\end{displaymath}
	where $S$ is uniformly persistent, i.e., there exists a constant $\epsilon$ where
	\begin{displaymath}
	\epsilon=\rho N\left(1-\frac{\frac{\alpha_{sa}}{\sigma}+\frac{\beta\rho N}{\sigma}}{\frac{\alpha_{sa}}{\sigma}+1+\frac{\beta\rho N}{\sigma}}\right)=\frac{\rho N}{\frac{\alpha_{sa}}{\sigma}+1+\frac{\beta\rho N}{\sigma}}
	\end{displaymath}such that 
	\begin{displaymath}
	\epsilon\leq\liminf\limits_{t\rightarrow\infty}S(t)\leq\limsup\limits_{t\rightarrow\infty}S(t)\leq \rho N.
	\end{displaymath}
	The persistence of $S$ leads to the persistence of $A$ and $L$. More specifically, \begin{displaymath}
	\liminf\limits_{t\rightarrow\infty}A(t) \geq \frac{\alpha_{sa}\epsilon}{(\alpha _{as}+\alpha_{al})  }=\epsilon_A \mbox{ and }\liminf\limits_{t\rightarrow\infty}L(t) \geq \frac{\alpha_{al}\epsilon_A}{(\alpha_{lc}+\alpha _{ls})}.
	\end{displaymath}
\end{theorem}
\textbf{Notes:} The technical proof of Theorem \ref{positiveinvariant}  is provided in the supplementary material file. This theorem indicates that Model \eqref{colonymigration} is biologically well-defined. Note that, within $\frac{1}{\sigma}$ minutes, an worker can independently discover the new site $\alpha_{sa}$ times and successfully recruit $\beta\rho N$ searchers, where $\frac{1}{\sigma}$ is the he maximum duration of ants in their population.  Theorem \ref{positiveinvariant} implies that, for a colony with $\rho N$ active workers, there are always at least $\epsilon$ searchers who are outside and search for a better home. The minimum scale of persistent searchers $\epsilon$ is increasing with respect to the maximum duration time $\frac{1}{\sigma}$, and is decreasing with respect to the discovery rate $\alpha_{sa}$ and maximum recruitment rate $\beta$. Note that  $S+A+L+C=\rho N$ and population $P$ does not depend on populations $A$ and $L$, these properties allow us to simplify Model \eqref{colonymigration} as follows
\begin{equation}\label{colonymigrationnew}
\begin{aligned}
\frac{dA}{dt}&=\left(\alpha_{sa}+\beta_{ls}L+\beta_{cs} C\right)\left(\rho N-A-L-C\right)-\alpha_{as}A-\alpha_{al}  A,\\
\frac{dL}{dt}&=\alpha_{al}  A-\alpha_{lc} Q_1 L-\alpha _{ls}L,\\
\frac{dC}{dt}&=\alpha_{lc} Q_1 L-\alpha _{cs}C,
\end{aligned}
\end{equation}
with
\begin{displaymath}
\left\{
\begin{array}{lr}
Q_1=0, \quad \text{if} \quad A+L+C<\Theta, \\
Q_1=1, \quad \text{if} \quad A+L+C>\Theta.\\
\end{array}
\right.
\end{displaymath}
System \eqref{colonymigrationnew} is a Filippov system \cite{filippov1988equations,meza2005threshold,da2006application,boukal1999lyapunov} which can be converted to a generalized form. Let $H(Z)=A+L+C-\Theta$ with vector $Z=(A, L, C)^{T}$, and
\begin{displaymath}
F_{S_{1}}(Z)=\begin{pmatrix}
\left(\alpha_{sa}+\beta_{ls}L+\beta_{cs} C\right)\left(\rho N-A-L-C\right)-(\alpha _{as}+\alpha_{al})  A\\
\alpha_{al}  A-\alpha _{ls}L\\ -\alpha _{cs}C
\end{pmatrix},
\end{displaymath}
\begin{displaymath}
F_{S_{2}}(Z)=\begin{pmatrix}
\left(\alpha_{sa}+\beta_{ls}L+\beta_{cs} C\right)\left(\rho N-A-L-C\right)-(\alpha _{as}+\alpha_{al})  A\\
\alpha_{al}  A-(\alpha_{lc}+\alpha _{ls})L\\
\alpha_{lc} L-\alpha _{cs}C
\end{pmatrix}.
\end{displaymath}
Then System \eqref{colonymigrationnew} can be rewritten as the following generalized Filippov system
\begin{equation}\label{filippov}
\dot{Z}=\left\{
\begin{array}{lr}
F_{S_{1}}(Z), \quad Z \in S_{1}, \\
F_{S_{2}}(Z), \quad Z \in S_{2},\\
\end{array}
\right.
\end{equation}
where $S_1=\left\{Z\in\Gamma\mid H(Z)<0\right\},$  $S_2=\left\{Z\in\Gamma\mid H(Z)>0\right\}$ are two regions divided by the discontinuity manifold
\begin{displaymath}
\Sigma=\left\{Z\in \Gamma \mid H(Z)=0\right\},
\end{displaymath}
and $\Gamma= \left\{(A, L, C)\mid 0\leq A+L+C\leq \rho N \right\}.$ We call System \eqref{filippov} defined in region $S_1$ as \textit{failed emigration state} and call System \eqref{filippov} defined in region $S_2$ as \textit{successful emigration state}.  The state portrait of System \eqref{filippov} is composed of the state portrait on $\Sigma$ and the state portraits in each regions $S_{i}$. Thus, we first study the dynamics of subsystems and the sliding mode on $\Sigma$ respectively.

\subsection{Dynamics of Subsystems and Equilibria of Filippov System \eqref{filippov}} Define 
$\eta_1=\frac{1}{\frac{\alpha_{sa}}{\alpha_{as}+\alpha_{al}}\left(1+\frac{\alpha_{al}}{\alpha_{ls}}\right)}, \quad
\xi_1=\frac{\frac{\alpha_{al}\beta_{ls}}{\alpha_{ls}(\alpha_{as}+\alpha_{al})}}{\frac{\alpha_{sa}}{\alpha_{as}+\alpha_{al}}\left(1+\frac{\alpha_{al}}{\alpha_{ls}}\right)},$ and 
\begin{displaymath}
L^{f}=\frac{\left(\rho N \xi_1-1-\eta_1\right)+\sqrt{\left(\rho N \xi_1-1-\eta_1\right)^2+4\xi_1\rho N}}{2\xi_1\left(1+\frac{\alpha_{ls}}{\alpha_{al}}\right)}, \quad A^{f}=\frac{\alpha_{ls}}{\alpha_{al}}L^{f}.
\end{displaymath}
Note that $\frac{\alpha_{sa}}{\alpha_{as}+\alpha_{al}}\left(1+\frac{\alpha_{al}}{\alpha_{ls}}\right)$ is the sum of times that a worker independently discovers new site and times that a worker transits from assessing population into leading population, within $\frac{1}{\alpha_{as}+\alpha_{al}}+\frac{1}{\alpha_{ls}}$ minutes, and $\frac{\alpha_{al}\beta_{ls}}{\alpha_{ls}(\alpha_{as}+\alpha_{al})}$ is times that a new leader recruits nestmates into new site in the same time period. Biologically, the interpretation to $\eta_1$ is the ‘recruitment efficiency’ of workers in new site without transportation recruitment, namely, during the average duration of workers, the ratio of the number of nestmates recruited by new leaders to the sum of numbers of new workers in each population (including assessing population and leading population). The interpretation to $\xi_1$ is the ‘input-output’ ratio of migrating colony, namely, during average duration of workers, the ratio of the initial number of searching workers to the sum of numbers of new workers in each population.\\

Let $\eta_2=\frac{1}{\frac{\alpha_{sa}}{\alpha_{as}+\alpha_{al}}\left[1+\frac{\alpha_{al}}{\alpha_{ls}+\alpha_{lc}}+\frac{\alpha_{al}\alpha_{lc}}{(\alpha_{ls}+\alpha_{lc})\alpha_{cs}}\right]}, \xi_2=\frac{\frac{\beta_{ls}\alpha_{al}}{(\alpha_{lc}+\alpha_{ls})}+\frac{\beta_{cs}\alpha_{al}\alpha_{lc}}{\alpha_{cs}(\alpha_{lc}+\alpha_{ls})}}{\alpha_{sa}\left[1+\frac{\alpha_{al}}{\alpha_{ls}+\alpha_{lc}}+\frac{\alpha_{al}\alpha_{lc}}{(\alpha_{ls}+\alpha_{lc})\alpha_{cs}}\right]},$ 	and 	
\begin{displaymath}
C^{s}=\frac{\rho N \xi_2-1-\eta_2+\sqrt{\left(\rho N \xi_2-1-\eta_2\right)^2+4\xi_2\rho N}}{2\xi_2\left[1+\frac{\alpha_{cs}}{\alpha_{lc}}+\frac{\alpha_{cs}\left(\alpha_{lc}+\alpha_{ls}\right)}{\alpha_{al}\alpha_{lc}}\right]}, L^{s}=\frac{\alpha_{cs}}{\alpha_{lc}}C^{s},
A^{s}=\frac{\alpha_{lc}+\alpha_{ls}}{\alpha_{al}}L^{s}.
\end{displaymath}
The interpretation to $\eta_2$ is the recruitment efficiency of workers in new site with transportation recruitment, namely, during average duration of workers, the ratio of the number of nestmates recruited by new leaders and carriers to the sum of numbers of new workers in each population (including assessing, leading and carrying population). $\xi_2$ also is an ‘input-output’ ratio of migrating colony, namely, during average duration of workers, the ratio of initial number of searching population to the sum of numbers of new workers in each population.
We have the following results regarding analyzing the Filippov system \eqref{filippov}:

\begin{theorem}\label{Lemmasubsystem12}
	If $A(t)+L(t)+C(t)<\Theta$, the Filippov system \eqref{filippov} becomes the  following model
	\begin{equation}\label{subsystemS1}
	\begin{aligned}
	\frac{dA}{dt}&=\left(\alpha_{sa}+\beta_{ls}L+\beta_{cs} C\right)\left(\rho N-A-L-C\right)-\alpha _{as}A-\alpha_{al}  A,\\
	\frac{dL}{dt}&=\alpha_{al}  A-\alpha _{ls}L,\\
	\frac{dC}{dt}&=-\alpha _{cs}C,
	\end{aligned}
	\end{equation}which has a unique boundary equilibrium $E^{f}(A^{f}, L^{f},0)$ that is globally asymptotically stable.
	If $A(t)+L(t)+C(t)>\Theta$, the Filippov system \eqref{filippov} becomes the  following model
	\begin{equation}\label{subsystemS2}
	\begin{aligned}
	\frac{dA}{dt}&=\left(\alpha_{sa}+\beta_{ls}L+\beta_{cs} C\right)\left(\rho N-A-L-C\right)-\alpha _{as}A-\alpha_{al}  A,\\
	\frac{dL}{dt}&=\alpha_{al}  A-\alpha_{lc} L-\alpha _{ls}L,\\
	\frac{dC}{dt}&=\alpha_{lc} L-\alpha _{cs}C.
	\end{aligned}
	\end{equation} which has a unique interior equilibrium $E^{s}(A^{s}, L^{s}, C^{s})$ that is locally asymptotically stable. Moreover, if $\alpha_{as}>\alpha_{ls}>\alpha_{cs}$ and $\beta_{cs}>\beta_{ls}$, then the equilibrium $E^{s}$ is globally asymptotically stable.
\end{theorem}
\noindent \textbf{Notes.} The technical proof of Theorem \ref{Lemmasubsystem12} is provided in the supplementary material file.   In the case of $A(t)+L(t)+C(t)>\Theta$, the steady state value of population $P$ is not unique, which is governed by the initial values of populations $C$ and $P$. Mathematically, this is an interesting result. However, in natural colonies, it is difficult to find carriers in new sites that have not been discovered by scouts.  Theorem \ref{Lemmasubsystem12} provides the local stability of interior equilibrium $E^{s}$ and the global stability of $E^{s}$ under sufficient conditions. Extensive numerical simulations suggest that the interior equilibrium $E^{s}$ is always globally asymptotically stable. Some typical simulations are shown in Figure SM7. Thus, we conjecture that $E^{s}$ is globally asymptotically stable for all parameters. However, it is difficult to testify this conjecture in theory due to the complexity of system.  Moreover, in this case System \eqref{colonymigration} has only one steady-state value $(1-\rho)N$ of passive population.\\

In order to proceed more dynamical results of our system, we provide some definitions related to equilibrium in piecewise smooth system \cite{di2008bifurcations,kuznetsov2003one} as follows:
\begin{definition}[Regular equilibrium]
	A point $Z^{*}$ is called a regular equilibrium of System \eqref{filippov} if $F_{S_1}(Z^{*})=0$, $H(Z^{*})<0$ or $F_{S_2}(Z^{*})=0$, $H(Z^{*})>0$.
\end{definition}
\begin{definition}[Virtual equilibrium]
	A point $Z^{*}$ is called a virtual equilibrium of System \eqref{filippov} if $F_{S_1}(Z^{*})=0$, $H(Z^{*})>0$ or $F_{S_2}(Z^{*})=0$, $H(Z^{*})<0$.
\end{definition}
Define \begin{equation}\label{RiPlus}
\mathcal{N}_{i}:=\frac{\xi_i\Theta^2+\Theta(1+\eta_i)}{\rho\left(1+\xi_i\Theta\right)}=\frac{\Theta}{\rho}+\frac{\Theta\eta_i}{\rho\left(1+\xi_i\Theta\right)}
\end{equation}
The biological implication of $\mathcal{N}_{i}$ is one critical size of colony, at which the number of active workers in this colony is $\Theta$ plus the sum of workers (including assessors and leaders) that can fully recruit $\Theta$ nestmates into the new site before leaving. From \eqref{RiPlus}, the size of  $ \mathcal{N}_{i}$ is increasing with respect to the threshold value $\Theta$  and ‘recruitment efficiency’ $\eta_i$ , and is decreasing with respect to ‘input-output’ ratio $\xi_i$ and active worker ratio $\rho$.  
Then, we have the following results of equilibria for System \eqref{filippov}:
\begin{theorem}\label{onesiteath02}
	If $A^{f}+L^{f}<\Theta$, then the system \eqref{filippov} has a regular equilibrium $E^{f}_{R}(A^{f}_{R}, L^{f}_{R},0)$ (\textit{failed emigration state}), and if $A^{f}+L^{f}>\Theta$, then the system \eqref{filippov} has a virtual equilibrium $E_{V}^{f}(A^{f}_{V}, L^{f}_{V},0)$.
\end{theorem}
\noindent\textbf{Notes.} Theorem \ref{onesiteath02} gives sufficient conditions for the existence of regular equilibrium $E^{f}_{R}$ located in region $S_1$, namely, 
\begin{equation}\label{condition1}
A+L{|}_{E_{R}^{f}}=\frac{\left(\rho N \xi_1-1-\eta_1\right)+\sqrt{\left(\rho N \xi_1-1-\eta_1\right)^2+4\xi_1\rho N}}{2\xi_1}<\Theta.
\end{equation}
which is equivalent to $N<\mathcal{N}_{1}$. This condition indicates that the colony size $N$ has great impact on the dynamics of System \eqref{filippov}, namely, if $N<\mathcal{N}_{1}$ then the colony is more likely to stabilize at \textit{failed emigration state} $E^{f}_{R}(A^{f}_{R}, L^{f}_{R},0)$.

\begin{theorem}\label{onesitecth01}
	If $A^{s}+L^{s}+C^{s}>\Theta$, then the system \eqref{filippov} has a regular equilibrium $E_{R}^{s}(A_{R}^{s}, L_{R}^{s}, C_{R}^{s})$ (\textit{successful emigration state}), and if $A^{s}+L^{s}+C^{s}<\Theta$, then the system \eqref{filippov} has a virtual equilibrium $E_{V}^{s}(A_{V}^{s}, L_{V}^{s}, C_{V}^{s})$.
\end{theorem}
\noindent\textbf{Notes.} Theorem \ref{onesitecth01} implies that System \eqref{filippov} has a regular equilibrium $E_{R}^{s}$ located in region $S_2$ if the parameters meet
\begin{equation}\label{condition2}
A+L+C{|}_{E_{R}^{s}}=	\frac{\left(\rho N \xi_2-1-\eta_2\right)+\sqrt{\left(\rho N \xi_2-1-\eta_2\right)^2+4\xi_2\rho N}}{2\xi_2}>\Theta.
\end{equation}which  is equivalent to $N>\mathcal{N}_{2}$. This condition indicates that if  $N>\mathcal{N}_{2}$ then the colony is more likely to stabilize at \textit{successful emigration state}  $E_{R}^{s}(A_{R}^{s}, L_{R}^{s}, C_{R}^{s})$ where the passive population ($P$) could be completely moved into new site.

\subsection{Dynamics on threshold manifold $\Sigma$}

In order to investigate the dynamics on the separating manifold $\Sigma$, we first determine the existence of \textit{crossing set} and  \textit{sliding set} on $\Sigma$ by using Filippov convex method \cite{filippov1988equations,da2006application,boukal1999lyapunov,tang2012sliding,xiao2013dynamics,tang2012piecewise}.

Let $\gamma(Z)=\langle H_z(Z), F_{S_{1}}(Z)\rangle\langle H_z(Z), F_{S_{2}}(Z)\rangle,$ where $\langle \cdot \rangle$ denotes the standard scalar product and $H_{Z}(Z)$ is the non-vanishing gradient of smooth function $H$ on $\Sigma$. Define the \textit{crossing set} $\Sigma_{C}\subset\Sigma$ as 
\begin{displaymath}
\Sigma_{C}=\left\{Z\in\Sigma\mid \gamma(Z)> 0\right\},
\end{displaymath} 
and the \textit{sliding set} $\Sigma_{S}\subset\Sigma$ as 
\begin{displaymath}
\Sigma_{S}=\left\{Z\in\Sigma\mid \gamma(Z)\leq 0\right\}, 
\end{displaymath}
where $\Sigma_{S}=\Sigma\setminus\Sigma_{C}$. For System \eqref{filippov}, it is easy to get that
\begin{displaymath}
\gamma(Z)=\left[(\alpha_{sa}+\beta_{ls}L+\beta_{cs}C)(\rho N-A-L-C)-\alpha_{as}A-\alpha_{ls}L-\alpha_{cs}C\right]^2>0
\end{displaymath} 
for all $Z\in\Gamma$. Therefore, we have  follows
\begin{lemma}\label{lemmaSet}
	For system \eqref{filippov}, we have $\Sigma_{C}=\Sigma \quad \text{and} \quad \Sigma_{S}=\emptyset.$ \end{lemma}
\noindent\textbf{Notes.} According to the definitions of crossing and sliding set, if $Z_0\in\Sigma_{C}$, then the two vectors $F_{S_1}(Z_0)$ and $F_{S_2}(Z_0)$ point to the same side of $\Sigma$ (See Figure SM8a), and if $Z_0\in\Sigma_{S}$, then the vectors $F_{S_1}(Z_0)$ and $F_{S_2}(Z_0)$ point to the both side of $\Sigma$ (See Figure SM8b) or tangent to $\Sigma$. It indicates that the trajectories reaching $\Sigma_{C}$ immediately cross from one side to another, and the trajectories reaching $\Sigma_{S}$ may slide along the sliding vector (see Figure SM8b) to an internal point or the boundary of $\Sigma_{S}$. The result $\Sigma_{C}=\Sigma$ suggests that System \eqref{filippov} is a non-sliding piecewise system, i.e., all trajectories in System \eqref{filippov} hitting the manifold $\Sigma$ would cross into the opposite region instead of sliding on $\Sigma$. This implies that if System \eqref{filippov} has multiple locally stable regular equilibria, then the system has multiple attractors; while if System \eqref{filippov} has multiple virtual equilibria, then the system would likely to have oscillating dynamics. In the next section, we will classify dynamics of System \eqref{filippov} in more details.

\section{Dynamical behaviors of Filippov system \eqref{filippov}}\label{Section4} 
In this section, we explore the global dynamics of System \eqref{filippov}. It follows from Theorem \cref{onesiteath02} and Theorem \ref{onesitecth01} that System \eqref{filippov} can have zero, one and two regular equilibria according to the relationship between $N$ and $\mathcal{N}_i$ ($i=1,2$). Thus, based on the relationship between $N$ and $\mathcal{N}_i$ ($i=1,2$), we classify the possible dynamics of system in four cases which are provided in the following four corollaries, respectively.

\begin{corollary}\label{corollary1}
	System \eqref{filippov} has local stability at $E_{R}^{f}(A_{R}^{f}, L_{R}^{f},0)$ (\textit{failed emigration state}) if $N<\min\{\mathcal{N}_1, \mathcal{N}_2\}$.
\end{corollary}
If the colony size $N$ is small, System \eqref{filippov} has one regular equilibrium $E^{f}_{R}$. In this case, the trajectories starting from region $S_1$ tend to $E_{R}^{f}$, and the trajectories starting from region $S_2$ also tend to $E_{R}^{f}$ after they cross the separating manifold. Time series and phase plots of System \eqref{filippov} shown in Figure SM9 suggest that $E_{R}^{f}$ (\textit{failed emigration state}) is the unique attractor in this case. Biologically, if the size of colony is less than the sum of quorum threshold $\Theta$ and the number of active workers necessary to fully recruit $\Theta$ nestmates into new site, then the colony stabilize at the failed emigration state.

\begin{corollary}\label{corollary2}
	System \eqref{filippov} has local  stability at $E_{R}^{s}(A_{R}^{s}, L_{R}^{s}, C_{R}^{s})$ (\textit{successful emigration state}) if $N>\max\{\mathcal{N}_1, \mathcal{N}_2\}$.
\end{corollary}

If the colony size $N$ is large enough, System \eqref{filippov} has one regular equilibrium $E^{s}_{R}$. In this case, all solutions of System \eqref{filippov} tend to the equilibrium $E_{R}^{s}$ as shown in Figure SM10. Biologically, if the size of colony is greater than the sum of quorum threshold $\Theta$ and the number of active workers necessary to fully recruit $\Theta$ nestmates into new site, then the colony reach consensus on emigration without splitting.

\begin{corollary}\label{corollary3}
	System \eqref{filippov} has only virtual equilibrium if $\mathcal{N}_1<N< \mathcal{N}_2$.
\end{corollary}

\begin{figure}[htbp]
	\centering
	\includegraphics[width=12.5cm]{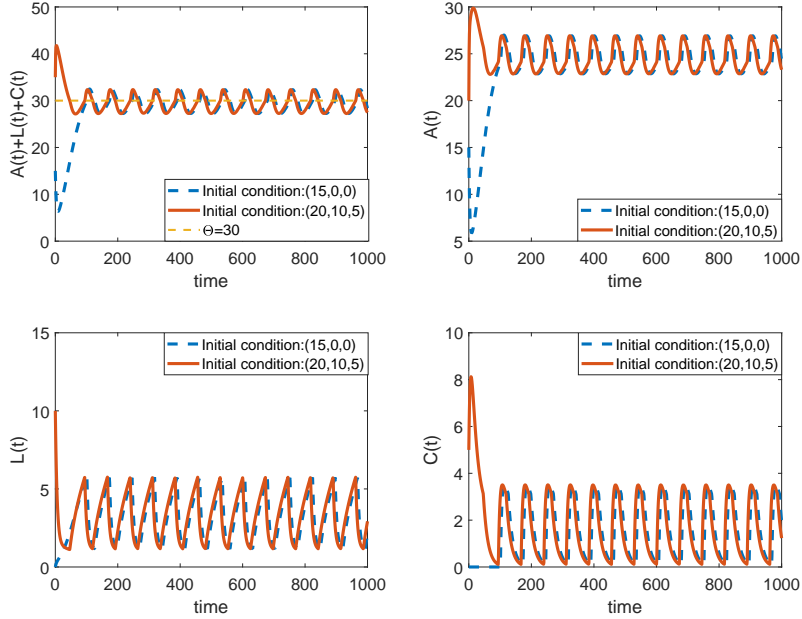}
	\caption{Time series plot shows the existence of oscillation when System \eqref{filippov} has no regular equilibrium. The parameters are  $N=200$, $\Theta=30$, $\rho=0.25$, $\alpha_{cs}=0.07$, $\alpha_{ls}=0.018$, $\beta_{ls}=0.049$, $\beta_{cs}=0.079$, $\alpha_{al}=0.007$, $\alpha_{lc}=0.15$, $\alpha_{as}=0.24$, $\alpha_{sa}=0.01$.}
	\label{TimeseriesScenarioC}
\end{figure}

From Corollary \ref{corollary3}, both $E^{f}$ and $E^{s}$ are virtual equilibria when the colony size $N$ is intermediate. Figure \ref{TimeseriesScenarioC} shows that, regardless of initial conditions, the size of total active workers at new site ($A(t)+L(t)+C(t)$) continuously oscillates around quorum threshold, and the oscillations are also found in each active population. Figure SM11 further illustrates that solutions starting from regions $S_1$ and $S_2$ tend to threshold interface, then go back and forth on both sides of the threshold interface along periodic orbits. In this case, System \eqref{filippov} constantly switches between failed emigration state and successful emigration state. Biologically, if the active workers in a colony with intermediate size can fully recruits $\Theta$ nestmates into new site before they leaving by using only tandem running , but cannot do so by using transportation, then the colony is undecided in the choice of new site and old nest.

\begin{corollary}\label{corollary4}
	System \eqref{filippov} has two regular equilibrium $E_R^{f} (A_{R}^{f}, L_{R}^{f},0)$ (\textit{failed emigration state}) and $E_R^{s} (A_{R}^{s}, L_{R}^{s}, C_{R}^{s})$ (\textit{successful emigration state}) which are always locally stable if $\mathcal{N}_2<N< \mathcal{N}_1$.
\end{corollary}

From Corollary \ref{corollary4}, in this case, both $E^{f}$ and $E^{s}$ are regular equilibria. It indicates that System \eqref{filippov} exhibits bistability between $E_{R}^{f}$ and $E_{R}^{s}$, namely the solutions with different initial conditions eventually stabilize at two levels (see Figure SM12). Biologically, if the active workers in a colony with intermediate size can fully recruits $\Theta$ nestmates into new site before they leaving by using transportation, but cannot do so by using only tandem running, then the colony either emigrate to the new finding site or maintain the original nest. We will analyze how the initial values affect the solutions of System \eqref{filippov} at bistable state in more details.\\

\begin{figure}[tbhp]
	\centering
	\subfloat[ ]{\label{Initialcondition:subfig1}\includegraphics[width=6cm]{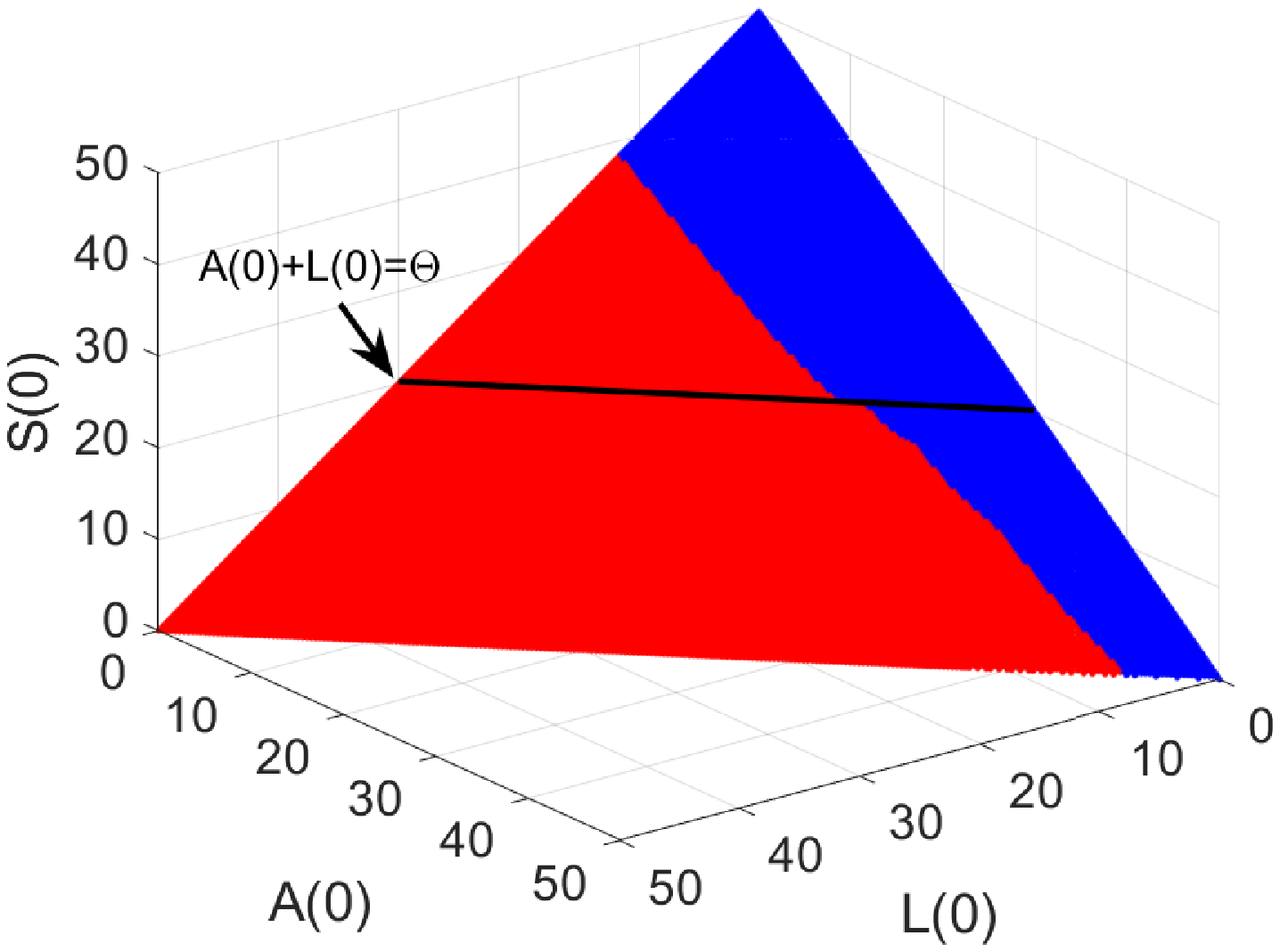}}\hspace{5mm}
	\subfloat[ ]{\label{Initialcondition:subfig2}\includegraphics[width=6cm]{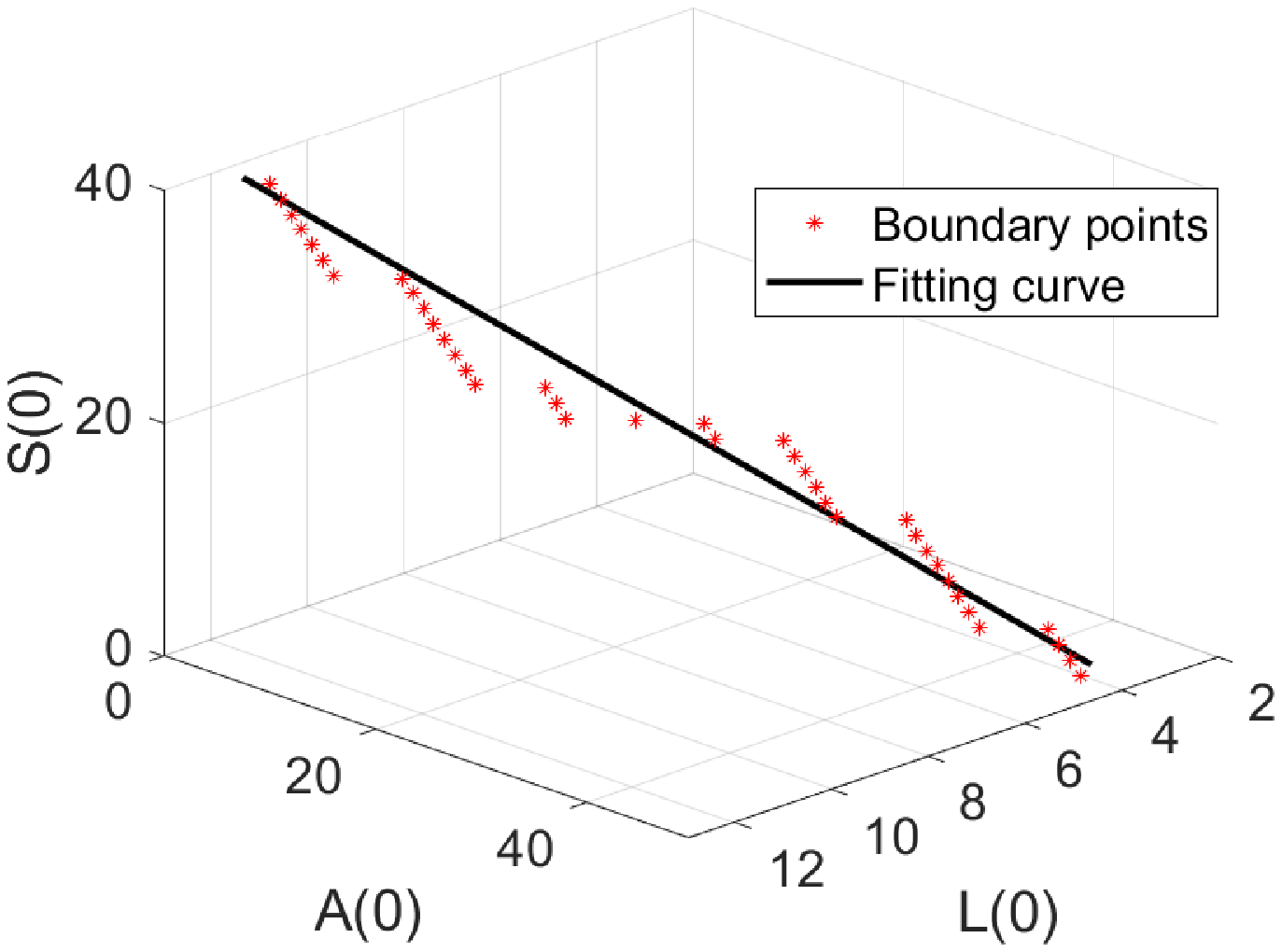}}
	\caption{Figure \ref{Initialcondition:subfig1} is the basin attractions of System \eqref{filippov} with the parameters taken as in Figure SM12 and $S(0), A(0), L(0)\in(0, \rho N)$, $C(0)=0$, where the red region is the basins of  atrractions of $E_{R}^{s}$, the blue region is the basins of atrractions of $E_{R}^{f}$. Figure \ref{Initialcondition:subfig2} is the fitting curve of boundary points btween two regions in Figure \ref{Initialcondition:subfig1}.}
	\label{Initialcondition}
\end{figure}

\noindent\textbf{Initial Condition Impact For Bistable Case:}
From Figure SM12, the solution starting from region $S_1$ reaches to equilibrium $E_R^{s}$ and the solution starting from region $S_2$ reaches to equilibrium $E_R^{f}$. It implies that, the relationship between initial values of active workers and quorum threshold does not completely determine whether the trajectory is tending to $E_R^{f}$ or $E_R^{s}$. In order to explore how do initial conditions affect the dynamics of System \eqref{filippov}, we take extensive numerical simulations to obtain an estimate of basins of attractions of System \eqref{filippov} with varying $S(0)$,  $A(0)$ and  $L(0)$ ($C(0)=0$). 
A typical simulation
is shown in  Figure \ref{Initialcondition:subfig1}. From Figure \ref{Initialcondition:subfig1}, we can obtain the following results: (1) If $L(0)=0$, all solutions tend to equilibrium $E_{R}^{f}$ regardless of the variations of $S(0)$ and $A(0)$; (2) If $L(0)>0$, the solutions with $S(0)+A(0)+L(0)<\Theta$ tend to equilibrium $E_{R}^s$ when $L(0)$ is large enough, and the solutions with $S(0)+A(0)+L(0)>\Theta$ tend to equilibrium $E_{R}^f$ when $L(0)$ is small.

In order to illustrate the importance of $L(0)$ on the outcomes of dynamics quantitatively, we fit the boundary between two basins of attractions of $E_R^{f}$ and $E_R^{s}$ as shown in Figure \ref{Initialcondition:subfig2}, where the red points are the boundary points on the basins of attractions of $E_R^{f}$ (red region in Figure \ref{Initialcondition:subfig1}) that connect with the basins of attractions of $E_R^{s}$ (blue region in Figure \ref{Initialcondition:subfig1}), the black line is the fitting curve of these red points. The function of the fitting curve is
\begin{equation}\label{fittingcurve1}
\frac{L(0)-a_1}{a_4}=\frac{A(0)-a_2}{a_5}=\frac{S(0)-a_3}{a_6},
\end{equation}
where $a_1=5.9739$, $a_2=31.0209$, $a_3=12.0327$, $a_4=0.8444$, $a_5=-4.8755$, $a_6=4.1945$.
It then follows from $\lvert a_4 \rvert \ll \lvert a_6 \rvert < \lvert a_5 \rvert$ that $L(0)$ has a much lower rate of change along the fitting curve than $S(0)$ or $A(0)$. This result indicates that, nearby the fitting curve, System \eqref{filippov} is more sensitive to the variations of $L(0)$ than the variations of $S(0)$ or $A(0)$. From $a_1<a_3<a_2$, the values of $L(0)$ are much less than the values of $S(0)$ or $A(0)$ along the fitting curve. It indicates that the solutions of System \eqref{filippov} with larger $L(0)$ are much more likely to tend to equilibrium $E_R^{s}$. We also take extensive numerical simulations of basins of attractions with varying $S(0)$, $A(0)$ and $C(0)$ ($L(0)=0$), as well as fit the boundary curves. The results indicate that the solutions with larger $C(0)$ are much more likely to tend to equilibrium $E_R^{s}$.

Figure \ref{Initialcondition} suggests that the initial values of recruiters (including leaders and carriers) have greatly impact on dynamical patterns when System \eqref{filippov} is in bistable state. From the biological point of view, if sudden environmental disturbance kills abundant active ants who are migrating from old nest to new finding site, the size of surviving recruiters at new site plays a crucial role in the decision to keep migrating. 

\section{Synergistic effects of colony size and quorum threshold on the dynamics}\label{Section5}
In this section, we will explore the synergistic effects of colony size $N$ and quorum threshold $\Theta$ on the dynamics of System \eqref{filippov} by analysis and bifurcation approaches. 
Denote a critical size of recruiters 
\begin{displaymath}
\Theta_c=\frac{\alpha_{sa}(1-\frac{\alpha_{cs}}{\alpha_{ls}})}{\beta_{ls}(\frac{\alpha_{cs}}{\alpha_{ls}}-\frac{\beta_{cs}}{\beta_{ls}})}
\end{displaymath}
whose existence requires the same sign of $1-\frac{\alpha_{cs}}{\alpha_{ls}}$ and $\frac{\alpha_{cs}}{\alpha_{ls}}-\frac{\beta_{cs}}{\beta_{ls}}$, i.e., $1>\frac{\alpha_{cs}}{\alpha_{ls}}>\frac{\beta_{cs}}{\beta_{ls}}$ or $1<\frac{\alpha_{cs}}{\alpha_{ls}}<\frac{\beta_{cs}}{\beta_{ls}}$. Biologically, the existence of positive $\Theta_c$ is determined by the transition rates of two types of recruiters (leader and carrier) to search group  $S$ and the recruitment rates of two types of recruiters (leader and carrier) from search group $S$ to assessor group $A$. 
Recall that 
\begin{displaymath}
\mathcal{N}_{i} :=\frac{\xi_i\Theta^2+\Theta(1+\eta_i)}{\rho\left(1+\xi_i\Theta\right)}=\frac{\Theta}{\rho}+\frac{\Theta\eta_i}{\rho\left(1+\xi_i\Theta\right)}
\end{displaymath} which is increasing in $\Theta$ and $\eta_i$, and is decreasing in $\rho$ and $\xi_i$. In the following, we show how the existence of positive $\Theta_c$ is related to the relationship of $\mathcal{N}_1$ and $\mathcal{N}_2$ as follows:
\begin{theorem}\label{lemma1} If $\Theta_c<0$, then we have the following two cases:
	\begin{itemize}
		\item [(a)] If $\frac{\alpha_{cs}}{\alpha_{ls}}<\min\{1,\frac{\beta_{cs}}{\beta_{ls}}\}$, then  $\mathcal{N}_1(\Theta)>\mathcal{N}_2(\Theta)$ for all $\Theta>0$;
		\item [(b)] If $\frac{\alpha_{cs}}{\alpha_{ls}}>\max\{1,\frac{\beta_{cs}}{\beta_{ls}}\}$, then $\mathcal{N}_1(\Theta)<\mathcal{N}_2(\Theta)$ for all $\Theta>0$;
	\end{itemize}
	And if  $\Theta_c>0$, then we have the following two cases:		
	\begin{itemize}
		\item [(c)] If $1<\frac{\alpha_{cs}}{\alpha_{ls}}<\frac{\beta_{cs}}{\beta_{ls}}$, then $\mathcal{N}_1(\Theta)<\mathcal{N}_2(\Theta)$ for all $0<\Theta<\Theta_c$ and $\mathcal{N}_1(\Theta)>\mathcal{N}_2(\Theta)$  for all $\Theta>\Theta_{c}$;
		\item [(d)] If $1>\frac{\alpha_{cs}}{\alpha_{ls}}>\frac{\beta_{cs}}{\beta_{ls}}$, then $\mathcal{N}_1(\Theta)>\mathcal{N}_2(\Theta)$ for all $0<\Theta<\Theta_{c}$ and $\mathcal{N}_1(\Theta)<\mathcal{N}_2(\Theta)$ for all $\Theta>\Theta_{c}$.
	\end{itemize}
\end{theorem} 
\textbf{Notes.} The technical proof of Theorem \ref{lemma1} is provided in the supplementary materials.
Theorem \ref{lemma1} gives the relationships between $\mathcal{N}_1$ and $\mathcal{N}_2$ with respect to $\Theta$ and $\Theta_c$ under four scenarios  that are determined by the signs of $1-\frac{\alpha_{cs}}{\alpha_{ls}}$ and $\frac{\alpha_{cs}}{\alpha_{ls}}-\frac{\beta_{cs}}{\beta_{ls}}$. These results suggest that it is important to distinguish two populations of recruiters, $L$ and $C$,  in modeling migration process. In other words, if we consider all recruiters as a group, we are not able to capture the interaction between different recruitment methods or explain the complex dynamic behavior that may occur. Moreover, the existence of positive $\Theta_c$ suggests the co-existence of undecided case and bistability between  $E_{R}^{f}$ (\textit{failed emigration state}) and  $E_{R}^{s}$ (\textit{successful emigration state}) in $N$ and $\Theta$ space which will be shown in more details in the following.
\begin{figure}[ht]
	\centering
	\subfloat[ ]{\label{bifurcationNThetaTwo:subfig1}\includegraphics[width=6cm]{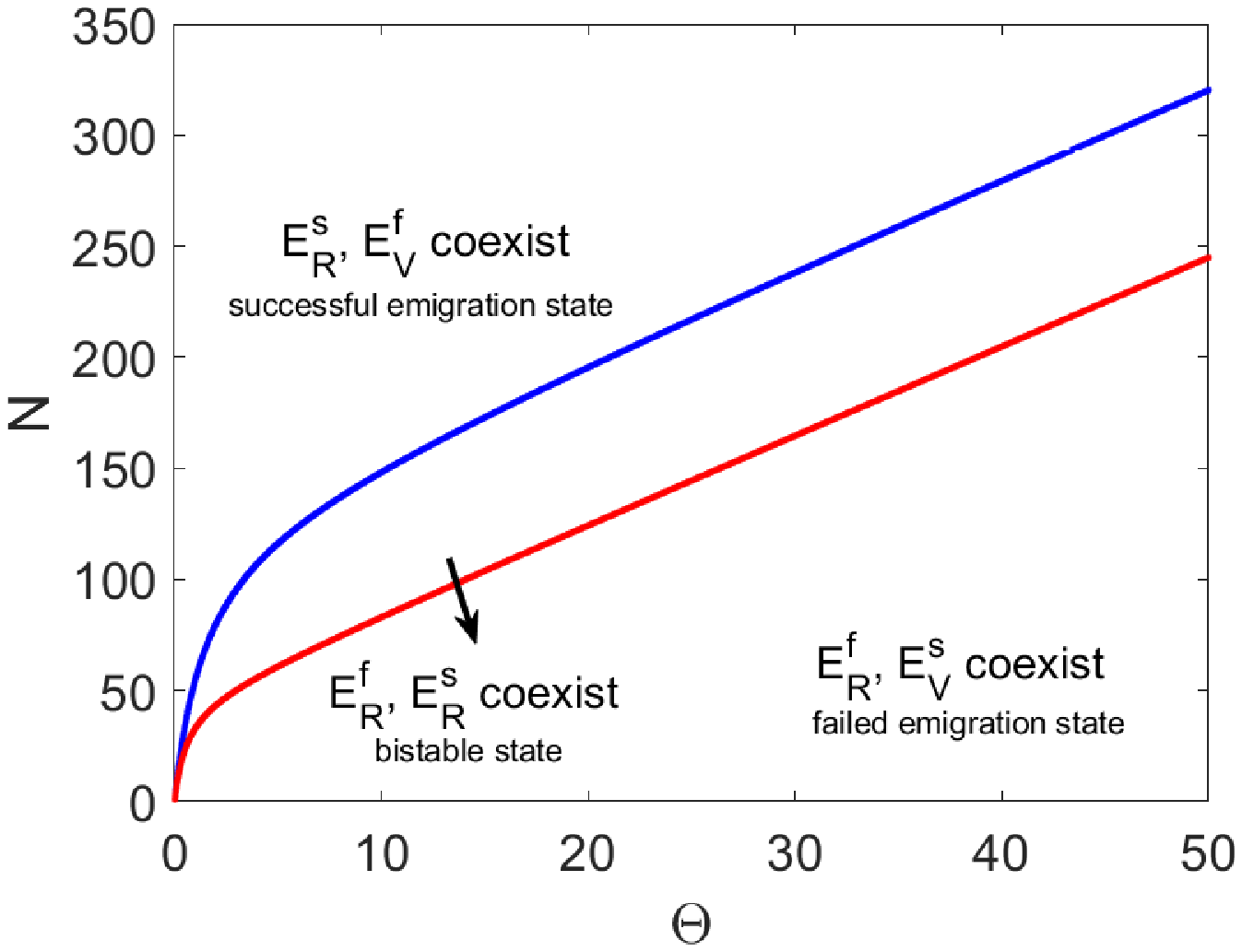}}
	\subfloat[ ]{\label{bifurcationNThetaTwo:subfig2}\includegraphics[width=6cm]{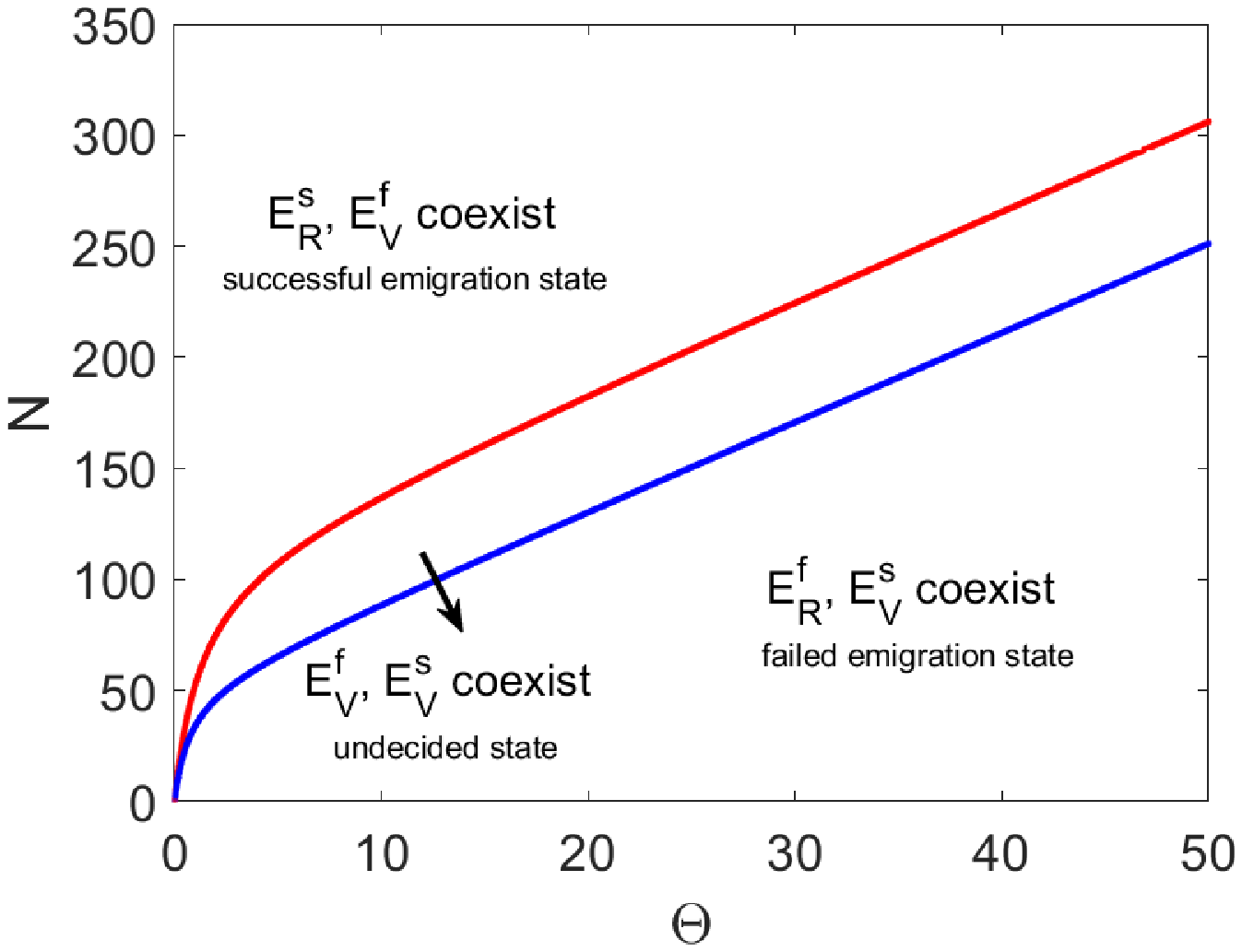}}\hspace{1mm}
	\subfloat[ ]{\label{bifurcationNThetaTwo:subfig3}\includegraphics[width=6cm]{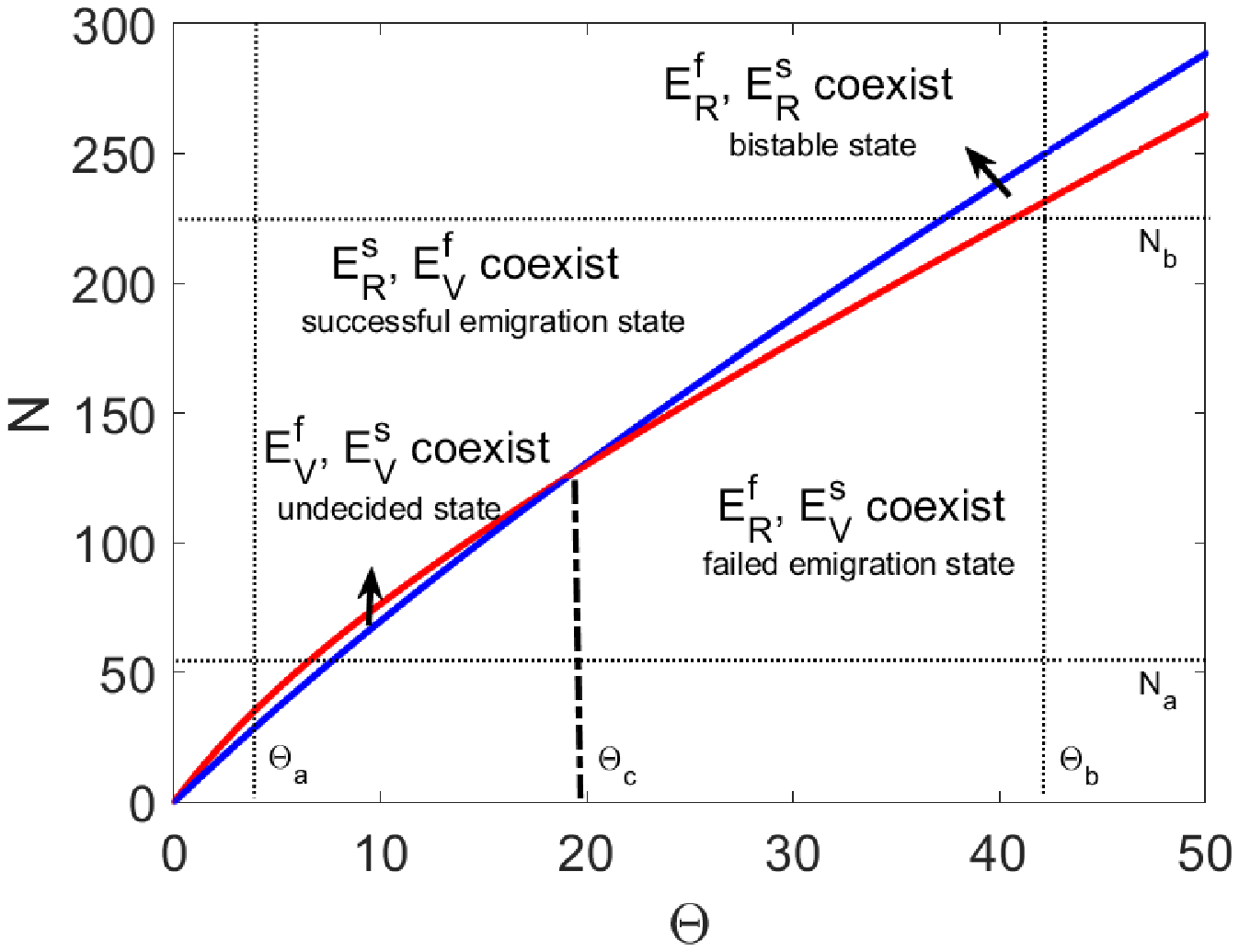}}
	\subfloat[ ]{\label{bifurcationNThetaTwo:subfig4}\includegraphics[width=6cm]{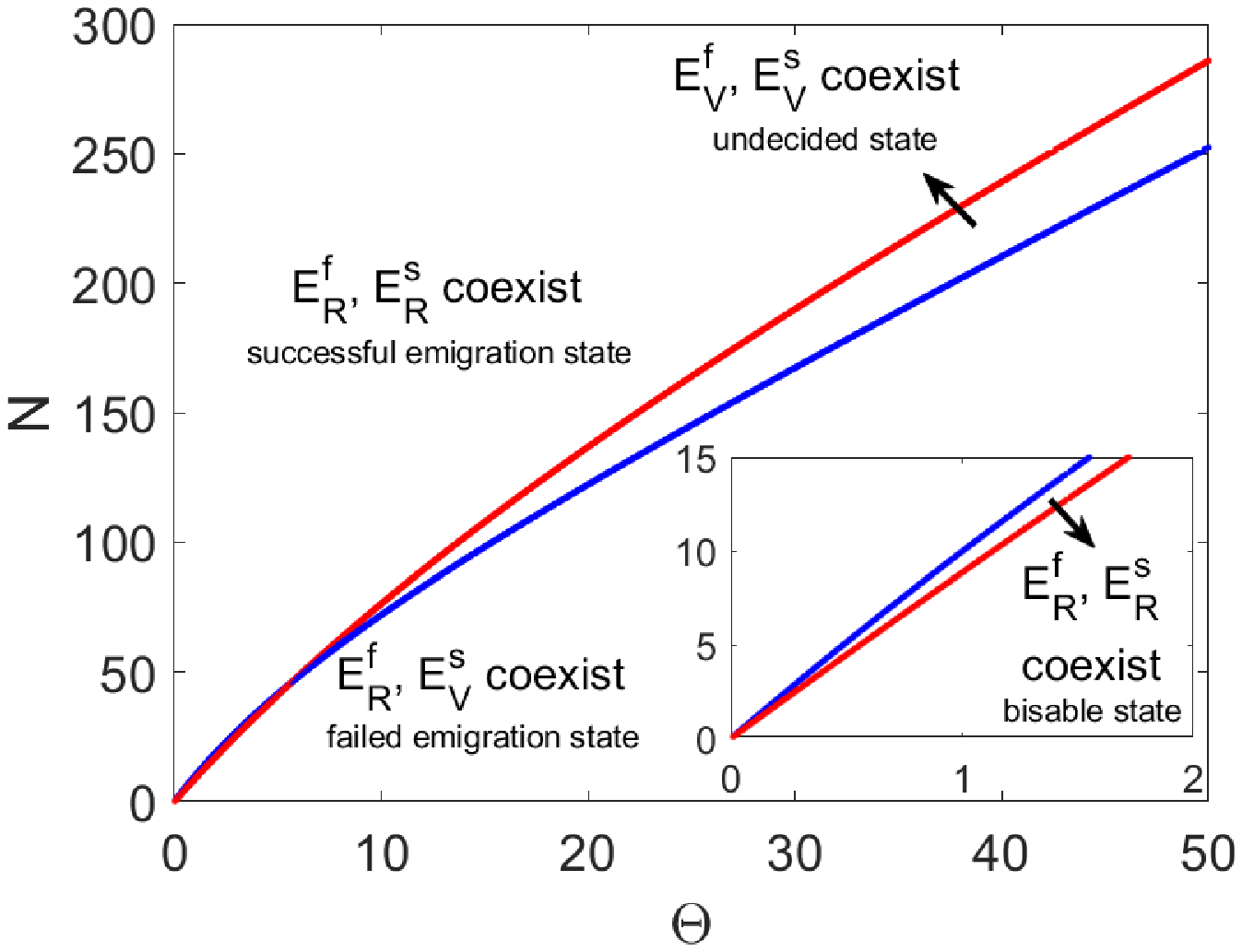}}
	\caption{Bifurcation diagrams of System \eqref{filippov} with respect to $N$ and $\Theta$ in four cases: (a) $\frac{\alpha_{cs}}{\alpha_{ls}}<\min\{1,\frac{\beta_{cs}}{\beta_{ls}}\}$; (b) $\frac{\alpha_{cs}}{\alpha_{ls}}>\max\{1,\frac{\beta_{cs}}{\beta_{ls}}\}$; (c) $1<\frac{\alpha_{cs}}{\alpha_{ls}}<\frac{\beta_{cs}}{\beta_{ls}}$; (d) $1>\frac{\alpha_{cs}}{\alpha_{ls}}>\frac{\beta_{cs}}{\beta_{ls}}$. Blue line is $\mathcal{N}_1(\Theta)$, red line is $\mathcal{N}_2(\Theta)$.
	}
	\label{bifurcationNThetaTwo}
\end{figure}

Based on Theorem \ref{lemma1} and the \cref{corollary1,corollary2,corollary3,corollary4}, we can obtain four possible regular/virtual equilibrium bifurcations of System \eqref{filippov} with respect to $N$ and $\Theta$.

\textit{Case (a)} $\frac{\alpha_{cs}}{\alpha_{ls}}<\min\{1,\frac{\beta_{cs}}{\beta_{ls}}\}$.

The $N$ and $\Theta$ parameter space is divided into three regions by curves $\mathcal{N}_1(\Theta)$ and $\mathcal{N}_2(\Theta)$. The existence of regular or virtual equilibrium in each region is indicated in Figure \ref{bifurcationNThetaTwo:subfig1}.  Figure \ref{bifurcationNThetaTwo:subfig1} suggests that System \eqref{filippov} always has at least one regular equilibrium, i.e., undecided state does not exist in this case.  

\textit{Case (b)} $\frac{\alpha_{cs}}{\alpha_{ls}}>\max\{1,\frac{\beta_{cs}}{\beta_{ls}}\}$.

The $N$ and $\Theta$ parameter space is also divided into three regions. The existence of equilibria in each region is indicated in Figure \ref{bifurcationNThetaTwo:subfig2}. From Figure \ref{bifurcationNThetaTwo:subfig2}, System \eqref{filippov} has at most one regular equilibrium, i.e., the bistability between $E_R^{f}$ and $E_{R}^s$ does not exist in this case.

\textit{Case (c)} $1<\frac{\alpha_{cs}}{\alpha_{ls}}<\frac{\beta_{cs}}{\beta_{ls}}$.

The $N$ and $\Theta$ parameter space is divided into four regions as shown in Figure \ref{bifurcationNThetaTwo:subfig3}. The existence of equilibria in each region implies that System \eqref{filippov}  has zero to two regular equilibria, i.e., System \eqref{filippov} has four possible dynamics (see the \cref{corollary1,corollary2,corollary3,corollary4}) in this case. Note that, if $\Theta<\Theta_c$, then undecided case exists for some $N$; and if $\Theta>\Theta_c$, then bistability exists for some $N$.

\textit{Case (d)} $1>\frac{\alpha_{cs}}{\alpha_{ls}}>\frac{\beta_{cs}}{\beta_{ls}}$.

The $N$ and $\Theta$ parameter space is also divided into four regions. The existence of equilibria shown in Figure \ref{bifurcationNThetaTwo:subfig4} is similar to case (c) but has an obvious difference, i.e., if $\Theta<\Theta_c$, bistability exists for some $N$, and if $\Theta>\Theta_c$, undecided case exists for some $N$.

\begin{figure}[ht]
	\centering
	\subfloat[ ]{\label{bifurcationNTheta:subfig1}\includegraphics[width=6cm]{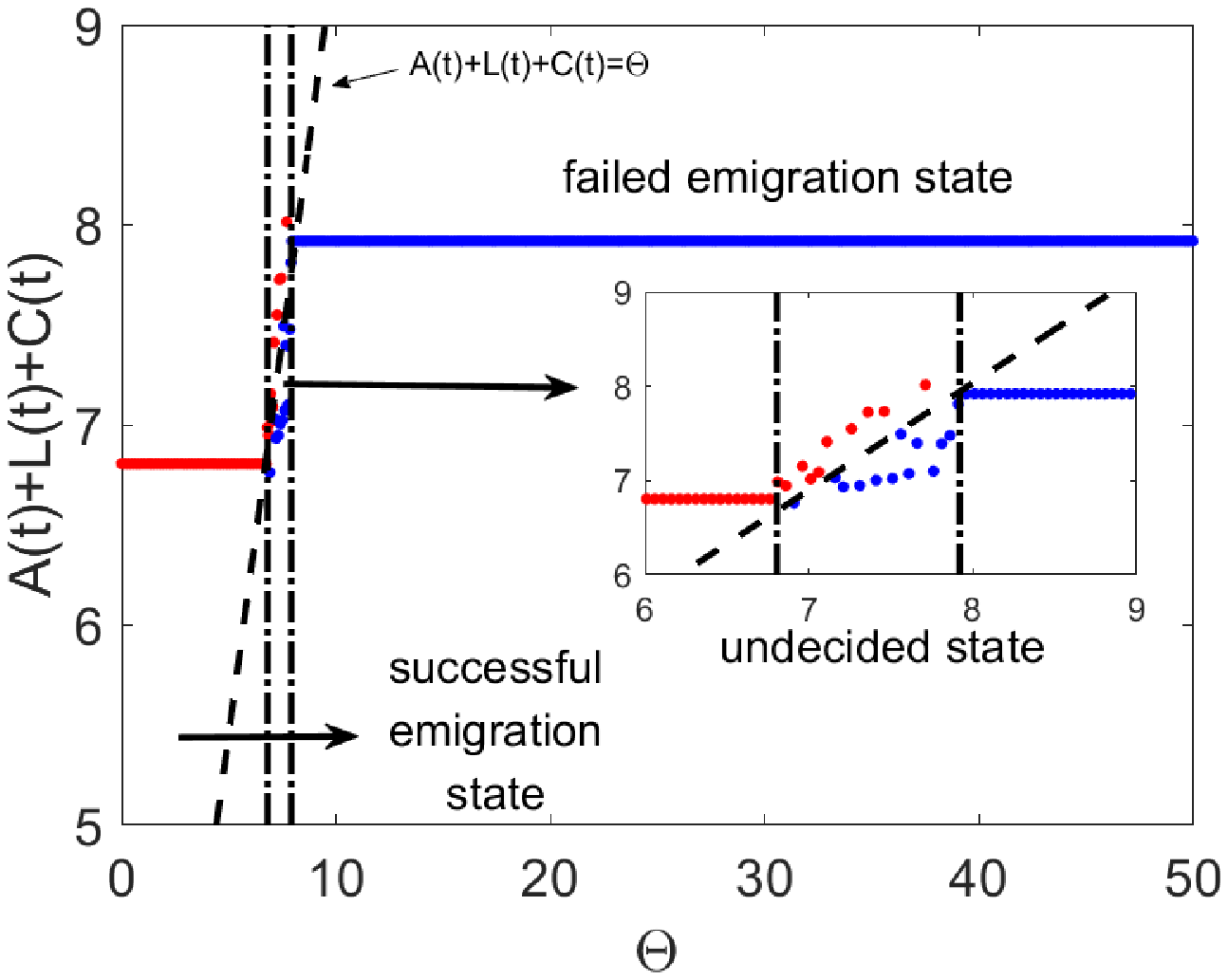}}
	\subfloat[ ]{\label{bifurcationNTheta:subfig2}\includegraphics[width=6cm]{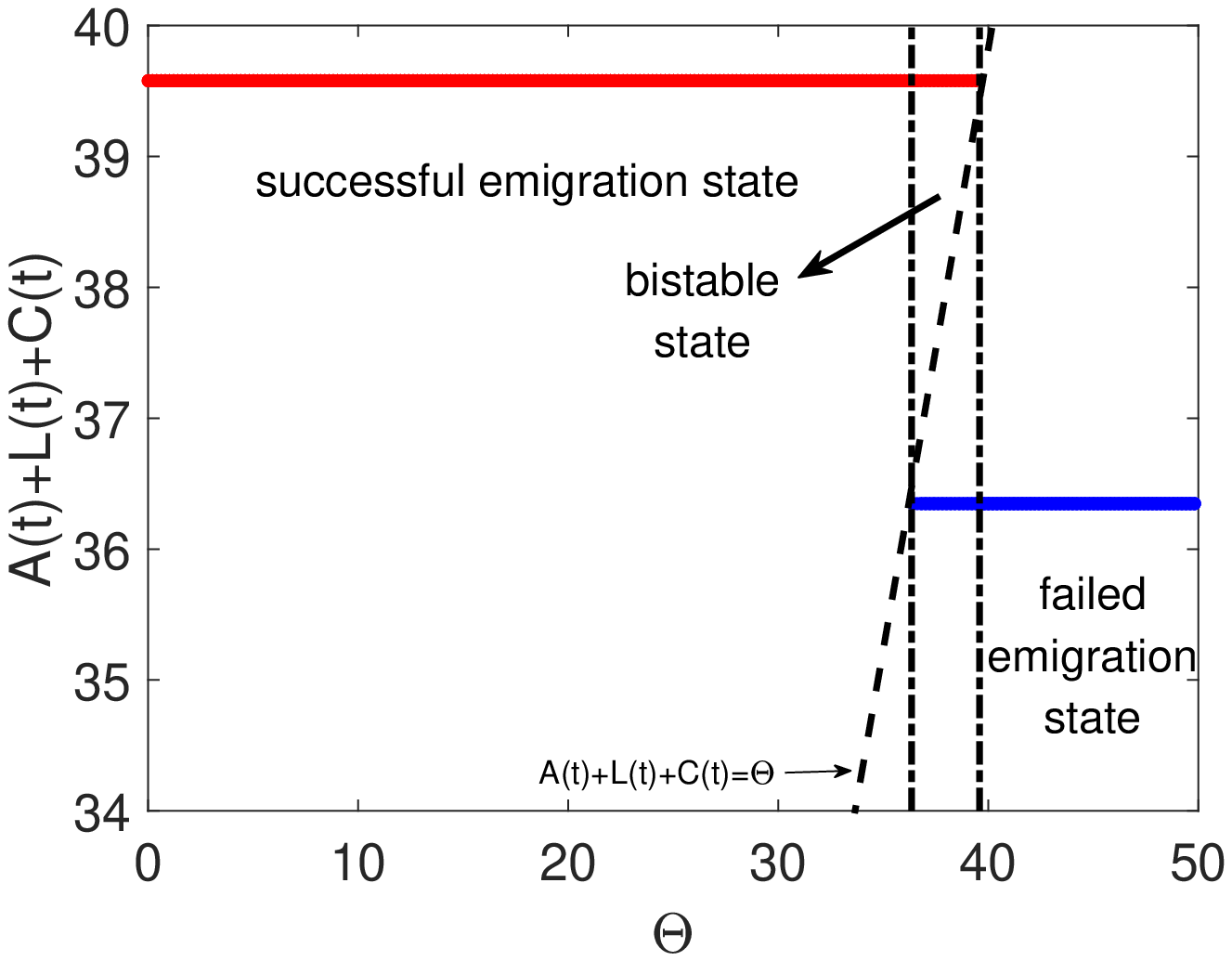}}\hspace{1mm}
	\subfloat[ ]{\label{bifurcationNTheta:subfig3}\includegraphics[width=6cm]{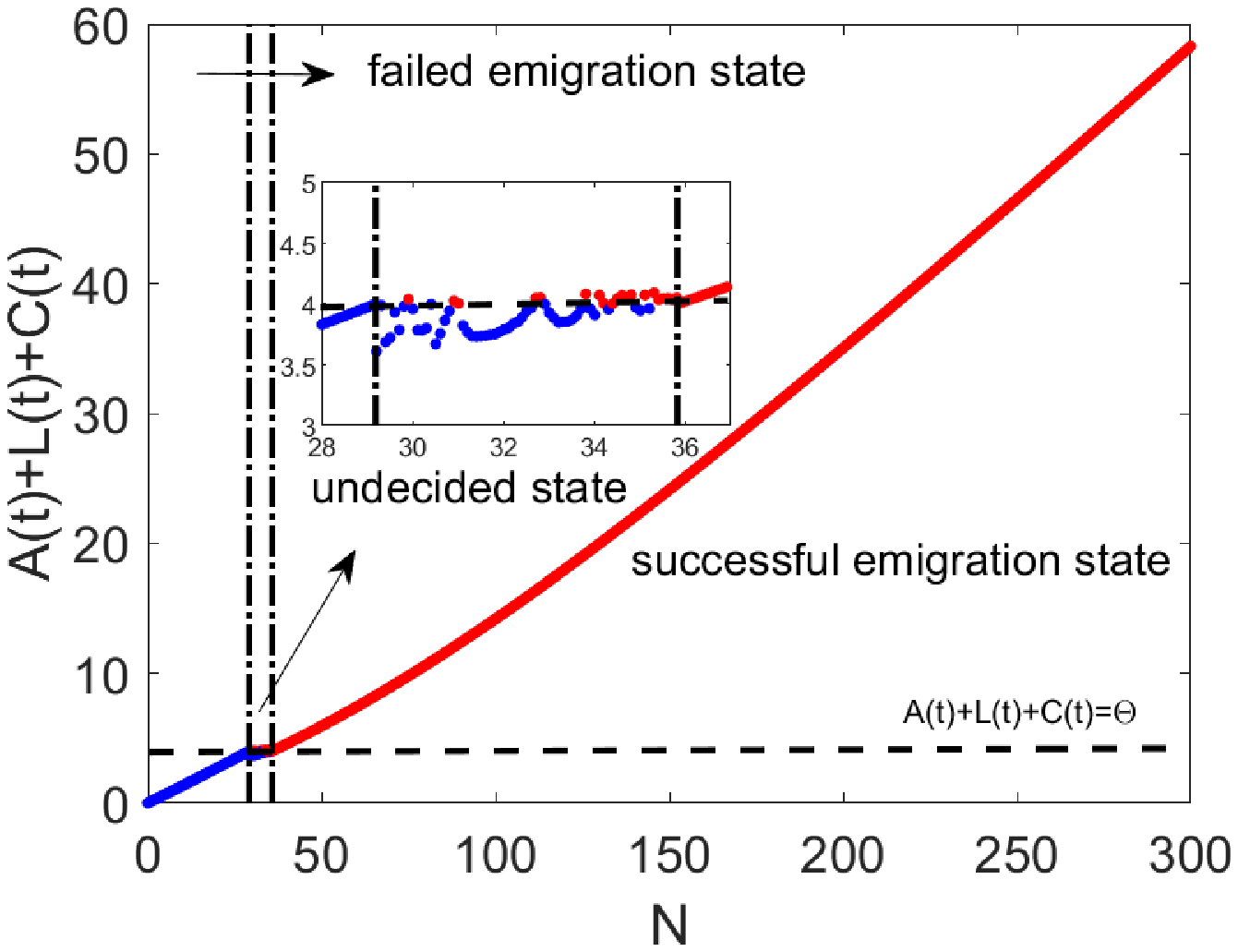}}
	\subfloat[ ]{\label{bifurcationNTheta:subfig4}\includegraphics[width=6cm]{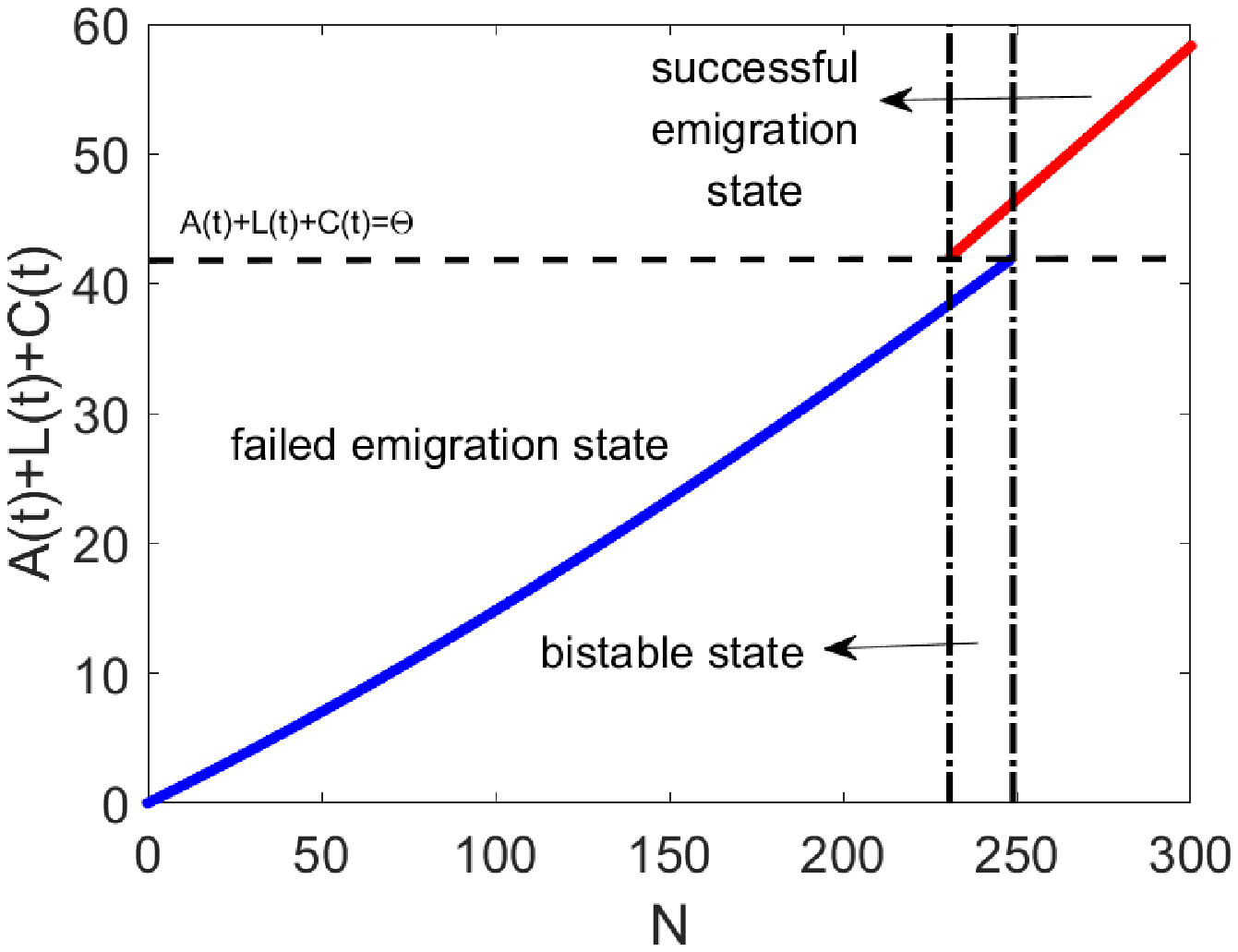}}
	\caption{Bifurcation diagrams of System \eqref{filippov} with two different levels of $N$ and two different levels of $\Theta$. In Figure \ref{bifurcationNTheta:subfig1}, $N=56$; In Figure \ref{bifurcationNTheta:subfig2}, $N=220$; In Figure \ref{bifurcationNTheta:subfig3}, $\Theta=4$; In Figure \ref{bifurcationNTheta:subfig4}, $\Theta=42$. The other parameters are $\rho=0.25$, $\alpha_{cs}=0.05$, $\alpha_{ls}=0.018$, $\beta_{ls}=0.004$, $\beta_{cs}=0.025$, $\alpha_{al}=0.057$, $\alpha_{lc}=0.28$, $\alpha_{as}=0.5$, $\alpha_{sa}=0.15$. }
	\label{bifurcationNTheta}
\end{figure}

In the following, we illustrate how does colony size and quorum threshold affect dynamics of System \eqref{filippov} in more details. We perform bifurcation study of System \eqref{filippov} satisfying $1<\frac{\alpha_{cs}}{\alpha_{ls}}<\frac{\beta_{cs}}{\beta_{ls}}$. We fix two different levels of $N$ (see $N_a$ and $N_b$ in Figure \ref{bifurcationNThetaTwo:subfig3}) and vary $\Theta$ to obtain bifurcation diagrams as shown in Figure \ref{bifurcationNTheta:subfig1} and Figure \ref{bifurcationNTheta:subfig2}, and fix two different levels of $\Theta$ (see $\Theta_a$ and $\Theta_b$ in Figure \ref{bifurcationNThetaTwo:subfig3}) and vary $N$ to obtain bifurcation diagrams as shown in  Figure \ref{bifurcationNTheta:subfig3} and Figure \ref{bifurcationNTheta:subfig4}. The bifurcation analysis of other cases can be obtained by the same way, which are provided detailed in supplementary materials.

For a colony with small level of size (see Figure \ref{bifurcationNTheta:subfig1}), when quorum threshold  is small (e.g., $\Theta$ varies from $0$ to $7$), System \eqref{filippov} stabilizes at $E_{R}^{s}$ (\textit{successful emigration state}); when quorum threshold is moderate (e.g., $\Theta$ varies from $7$ to $8$), points with two colors distributes discretely near the quorum threshold $\Theta$, i.e.,  System \eqref{filippov} is in undecided state; when quorum threshold is large (e.g., $\Theta$ varies from $8$ to $50$), System \eqref{filippov} stabilizes at $E_{R}^{f}$ (\textit{failed emigration state}). For a colony with large level of size (see Figure \ref{bifurcationNTheta:subfig2}), as $\Theta$ increases, the steady-state of colony also undergoes from successful emigration state (e.g., $\Theta$ varies from $0$ to $36$) to failed emigration state (e.g., $\Theta$ varies from $39$ to $59$) but with bistability between $E_{R}^{f}$ and $E_{R}^{s}$ as an intermediate (e.g., $\Theta$ varies from $36$ to $39$).

For small quorum threshold (see Figure \ref{bifurcationNTheta:subfig3}), when the colony size is small (e.g., $N$ varies from $0$ to $28$), System \eqref{filippov} stabilizes at $E_{R}^{f}$ (\textit{failed emigration state}); when the colony size is moderate (e.g., $N$ varies from $28$ to $36$), System \eqref{filippov} stabilizes is in undecided state; when the colony size is large (e.g., $N$ varies from $36$ to $300$), System \eqref{filippov} stabilizes at $E_{R}^{s}$ (\textit{successful emigration state}). For large quorum threshold (see Figure \ref{bifurcationNTheta:subfig4}), as $N$ increases, the steady-state of colony also undergoes from failed emigration state (e.g., $N$ varies from $0$ to $230$) to successful emigration state (e.g., $N$ varies from $250$ to $300$) but with bistability between $E_{R}^{f}$ and $E_{R}^{s}$ as an intermediate (e.g., $N$ varies from $230$ to $250$). 

\section{Conclusion}\label{Section6}
Social insects are considered as one of the evolutionarily most successful organisms on earth, which exhibit diverse decentralized organizations resulting from interactions among individuals and environment. Colony migration is a perfect example of collective decision-making, which causes great concern for entomologists and conservationists \cite{holldobler1990ants,visscher2007group}. Many studies have explored the decision rules and communication signals guiding the individual behaviors during colony migration  \cite{todd2000precis,dornhaus2004ants}, however, the underlying mechanisms at group level are less well understood. The observation of colony migration in previous study predicts that large colony size is necessary for the collective decision making, and the quorum threshold is not always correlated with group size.  How does the colony size affect outcomes of migration? How do synergies of colony size and quorum threshold regulate migration dynamic behaviors? To address these questions, we develop a piecewise system with a switching threshold and analyze the impact of key parameters (colony size and quorum threshold) on the dynamical patterns.

The dynamical features of our model are summarized as follows: In the absence and presence of transportation (see Theorem \ref{Lemmasubsystem12}), the colony migration systems both have only equilibrium dynamics (with a unique stable equilibrium respectively). However, the equilibrium dynamics of migration system with recruitment switching (see Theorem \ref{onesiteath02} and Theorem \ref{onesitecth01}) is more complicated. The system may admit regular/virtual equilibrium $E^{f}$ and regular/virtual equilibrium $E^{s}$ based on the relationship of the colony size $N$ and a critical size $\mathcal{N}_{i}$ ($i=1,2$) of this colony.

Mathematical results (see the \cref{corollary1,corollary2,corollary3,corollary4}) suggests how the colony size affects outcomes of migration. If the colony size is very small (i.e., $N<\min\{\mathcal{N}_{1}, \mathcal{N}_{2}\}$), the system would like to be stabilize at failed emigration state; If  the colony size is large enough (i.e., $N>\max\{\mathcal{N}_{1}, \mathcal{N}_{2}\}$), the system would like to be stabilize at successful emigration state; If the colony size is at critical level (i.e., $\mathcal{N}_{2}<N<\mathcal{N}_{1}$ or $\mathcal{N}_{1}<N<\mathcal{N}_{2}$), the system would like to be in undecided state or bistability between failed emigration state and successful emigration state. The undecided state is one of the interesting findings of our work (see Figure \ref{TimeseriesScenarioC} and Figure SM11), that is, the number of active workers presented at new site fluctuates around quorum threshold over time. It indicates that the colony switches between two sites and can not reach a consensus on nest selection. In fact, empirical studies \cite{mallon2001individual} has shown that the workers hesitate to recruit to poor sites. Our work also shows that the initial value of recruiter (who recruits nestmates through tandem running or transportation) plays an important role in determining which state the colony eventually tends to when system exhibits bistability (see Figure \ref{Initialcondition}). This result provides support to previous experimental studies \cite{pratt2002quorum} showing that tandem running and transportation recruitment offers great advantages for efficient emigration. More over, from the view on competition, System \eqref{filippov} can also be interpreted as the competition between old nest and new site for colonies. Specially, four dynamical patterns of System \eqref{filippov} have following explanations: (a) new site wins; (b) old nest wins; (c) no winner; (d) both sites have a potential to win. It provides a great new sight into understanding the decision-making issues on colony migration in social insects.

Bifurcation analysis (see Figure \ref{bifurcationNThetaTwo} and Figure \ref{bifurcationNTheta}) reveals how the synergies of colony size and quorum threshold regulate the dynamics of migration system \eqref{filippov}. If the quorum threshold is relatively low to colony size, then System \eqref{filippov} is more likely to stabilize at \textit{successful emigration state}. If the quorum threshold is relatively high to colony size, then System \eqref{filippov} is able to stabilize at \textit{failed emigration state}. The dynamics of System \eqref{filippov} with relative intermediate quorum threshold is more complicated, which is also determined by the critical size of recruiters ($\Theta_c$) as well as the recruitment rates and transition rates of two recruiters (signs of $1-\frac{\alpha_{cs}}{\alpha_{ls}}$ and $\frac{\alpha_{cs}}{\alpha_{ls}}-\frac{\beta_{cs}}{\beta_{ls}}$). Specially, if $\Theta_c<0$, then System \eqref{filippov} is in either undecided state or bistability between \textit{successful emigration state} and \textit{failed emigration state}, depending on the recruitment and transition rates of two recruiters; if $\Theta_c>0$, large colony and small colony would like to be in undecided state and bistable state respectively depending on the recruitment and transition rates. Our finding shows that the variations of colony size and quorum threshold greatly impact on migration. Empirical studies have claimed that the social insects could respond to environmental conditions or the need for urgency through adjusting their quorum \cite{mallon2001individual,pratt2002quorum}. For instance, the colonies will use a high quorum threshold to ensure a non-emergency and worthwhile emigration if the old nest remains intact, by contrast, they use a very small quorum threshold if the old nest is in a harsh situation \cite{dornhaus2004ants,franks2009speed}. In addition, our results may benefit experts interested in the potential factors influencing colony migration, such as transition rates and recruitment rates of two different recruiters.

In our current model, we assume that the quorum threshold is constant. This simplification allows us to obtain rigorous results on how colony size and quorum threshold affect the colony dynamics. However, this limitation also implies that our current model may not be a good description of the case that the quorum threshold could be correlated with colony size. Dornhaus et al. \cite{dornhaus2006colony} have shown that ants may measure the relative quorum, i.e., population in the new nest relative to that of the old nest, rather than the absolute number. Therefore, it is important to expand the colony migration model adopted relatively quorum threshold. The colony migration model is our first attempt. In addition to above suggestion, there are more reasonable and practical ways to extend this work. For instance: (i) In dynamical environment, the organisms are inevitably affected by environmental noise and demographic noise. It has been shown that the noises affect the interaction rate among group members and the follower's behavior in social insects \cite{al2013honey}. Thus, it would be interesting to incorporate the effect of randomly fluctuating environment in our model; (ii) \textit{Temnothorax} colonies can change the quorum size according to their colony size. They can achieve this end by considering the encounter rate at the old nest and at the target site. Thus, it would be a interesting subject to extend this model and investigate how the encounter rate affects collective decision making; (iii) In nature, the migrating social insects can evaluate several potential sites, compare them, and choose the best one, even most of the scouts visit only one site \cite{pratt2002quorum,pratt2005agent}. Therefore, it is interesting to propose a colony migration model with two or several potential sites to investigate the dynamic mechanism underlying nest-selection behavior, the effects of distances  or qualities on the outcome of migration, and the impact of colony size on the duration of collective decision-making. We keep these consideration for our future work.

\section*{Acknowledgments}
We have no conflict of interest to declare.


\bibliographystyle{unsrt}
\bibliography{references}

\end{document}


		
		
		
		

\begin{center}
	\normalfont\bfseries
	Supplementary Material for \textit{A Colony Migration-Decision Making Model With Hill Functions in Recruitment}
\end{center}


	
\section{Bifurcation analysis when $\frac{\alpha_{cs}}{\alpha_{ls}}<\min\{1,\frac{\beta_{cs}}{\beta_{ls}}\}$}\label{statementA}
	
\begin{figure}[htbp]
		\centering
		\includegraphics[width=6.5cm]{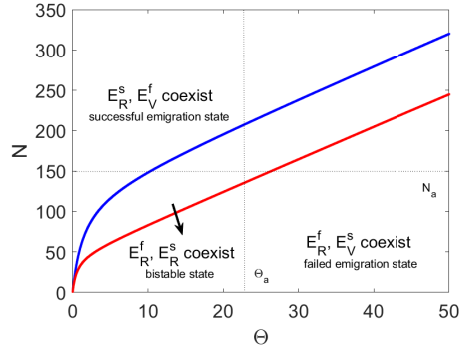}
		\caption{Bifurcation diagrams of System (3.2) with varying $N$ and $\Theta$ when $\frac{\alpha_{cs}}{\alpha_{ls}}<\min\{1,\frac{\beta_{cs}}{\beta_{ls}}\}$. Blue line is $\mathcal{N}_1(\Theta)$, red line is $\mathcal{N}_2(\Theta)$.}
		\label{Twoparametersbifurcation2}
	\end{figure}

\begin{figure}[htbp]
		\centering
		\subfloat[ ]{\label{bifurcationA:subfig1}\includegraphics[width=6cm]{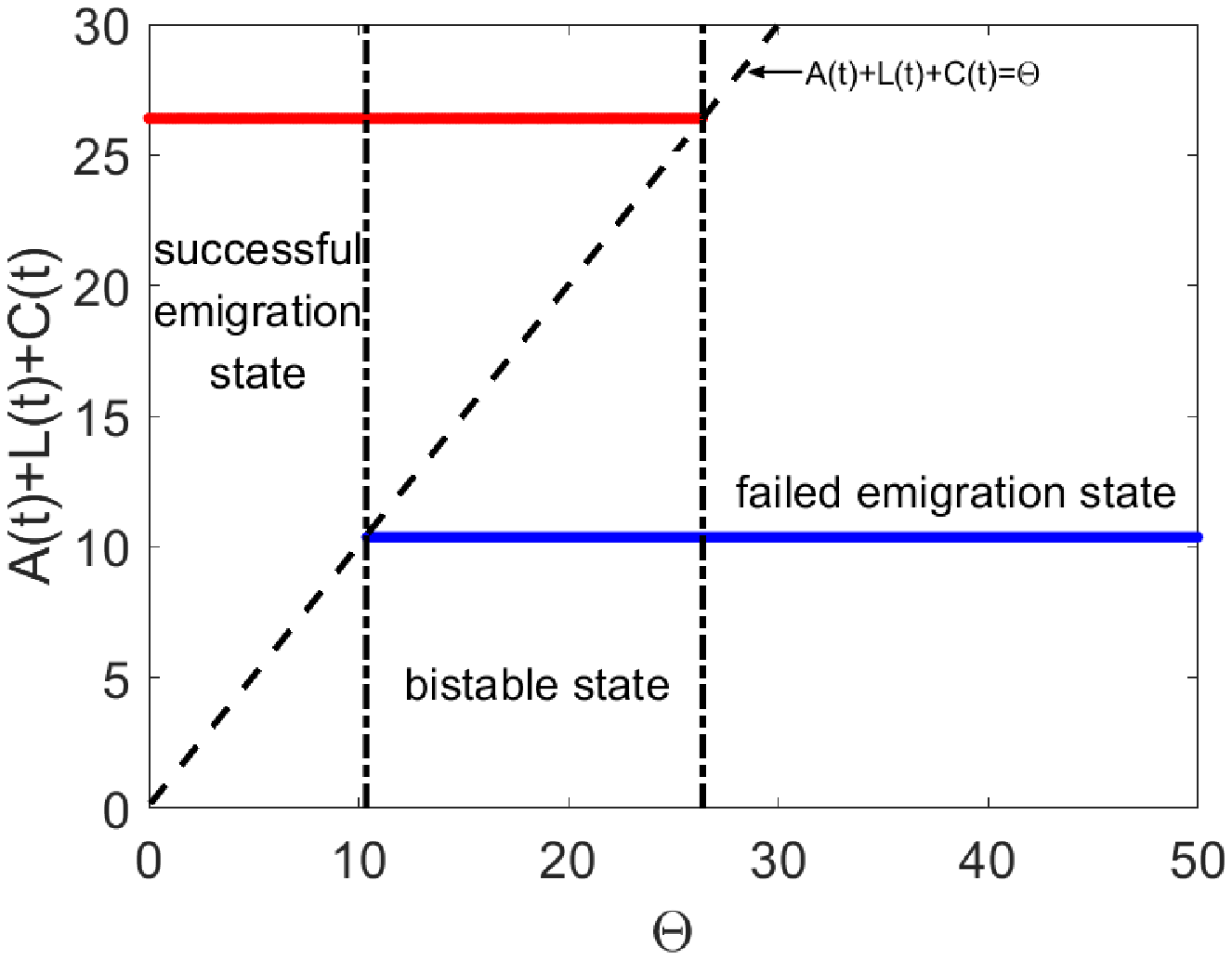}}
		\subfloat[ ]{\label{bifurcationA:subfig2}\includegraphics[width=6cm]{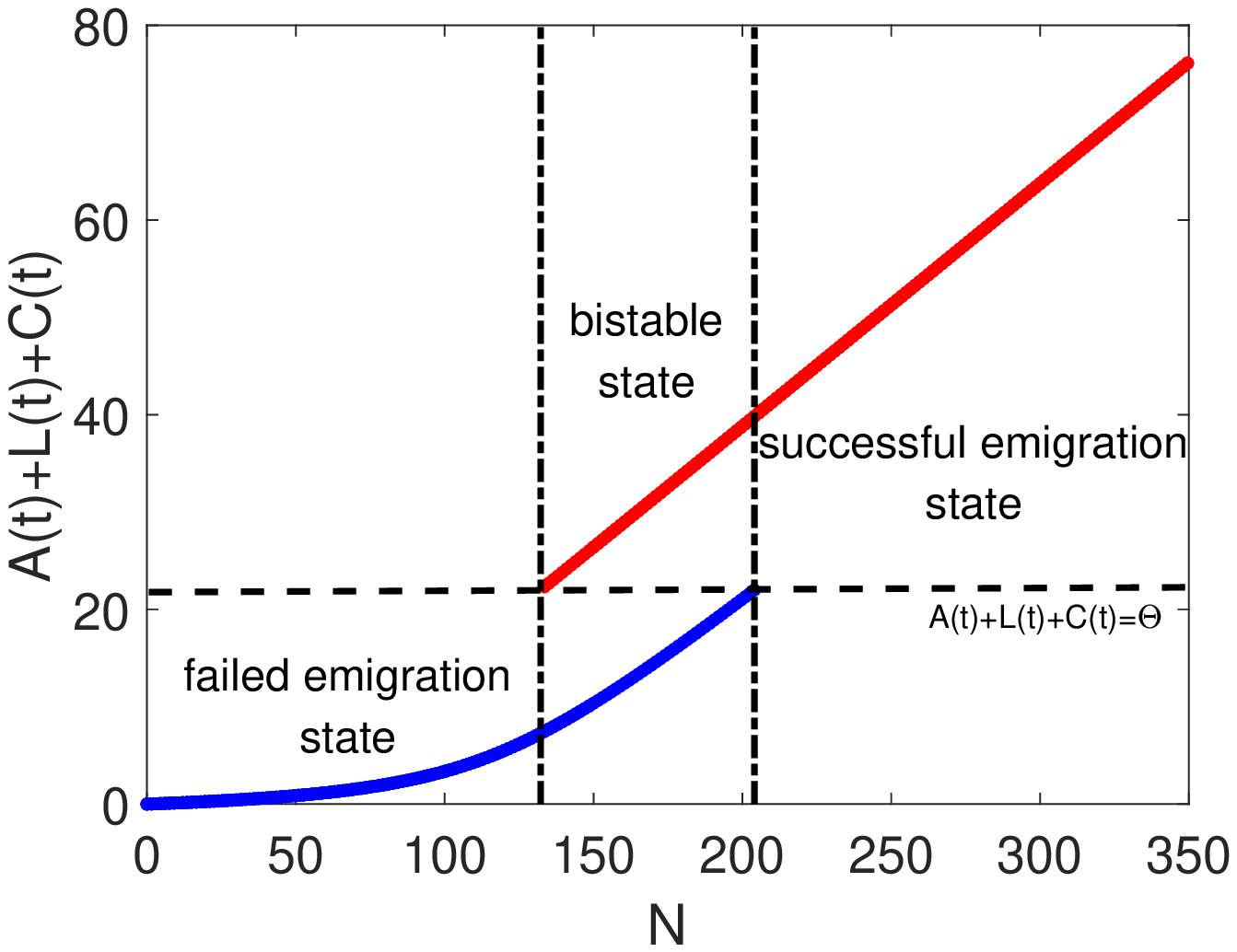}}
		\caption{Bifurcation diagrams of System (3.2). In Figure\ref{bifurcationA:subfig1}, $N=150$; In Figure \ref{bifurcationA:subfig2}, $\Theta=22$. The other parameters are $\rho=0.25$, $\alpha_{cs}=0.07$, $\alpha_{ls}=0.12$, $\beta_{ls}=0.033$, $\beta_{cs}=0.079$, $\alpha_{al}=0.032$, $\alpha_{lc}=0.15$, $\alpha_{as}=0.24$, $\alpha_{sa}=0.01$. }
		\label{bifurcationA}
	\end{figure}

	We perform bifurcation study of System (3.2) satisfying $\frac{\alpha_{cs}}{\alpha_{ls}}<\min\{1,\frac{\beta_{cs}}{\beta_{ls}}\}$. We fix $N=N_a$ in Figure \ref{Twoparametersbifurcation2} and vary $\Theta$ to obtain bifurcation diagrams as shown in Figure \ref{bifurcationA:subfig1}, and fix $\Theta=\Theta_a$ in Figure \ref{Twoparametersbifurcation2} and vary $N$ to obtain bifurcation diagrams as shown in Figure \ref{bifurcationA:subfig2}.
	
	For a colony with fixed colony size (see Figure \ref{bifurcationA:subfig1}	), when quorum threshold  is small (e.g., $\Theta$ varies from $0$ to $10$), System (3.2) stabilizes at $E_{R}^{s}$ (\textit{successful emigration state}); when quorum threshold is moderate (e.g., $\Theta$ varies from $10$ to $27$), System (3.2) is in bistability between $E_{R}^{f}$ and $E_{R}^{s}$; when quorum threshold is large (e.g., $\Theta$ varies from $27$ to $50$), System (3.2) stabilizes at $E_{R}^{f}$ (\textit{failed emigration state}). 
	
	For a fixed quorum threshold (see Figure \ref{bifurcationA:subfig2}), when colony size is small (e.g., $N$ varies from $0$ to $135$), System (3.2) stabilizes at $E_{R}^{f}$ (\textit{failed emigration state}); when colony size is moderate (e.g., $N$ varies from $135$ to $205$), System (3.2) is in bistability between $E_{R}^{f}$ and $E_{R}^{s}$; when colony size is large (e.g., $N$ varies from $205$ to $350$), System (3.2) stabilizes at $E_{R}^{s}$ (\textit{successful emigration state}). 
	
	\section{Bifurcation analysis when $\frac{\alpha_{cs}}{\alpha_{ls}}>\max\{1,\frac{\beta_{cs}}{\beta_{ls}}\}$}\label{statementB}

\begin{figure}[htbp]
		\centering
		\includegraphics[width=6.5cm]{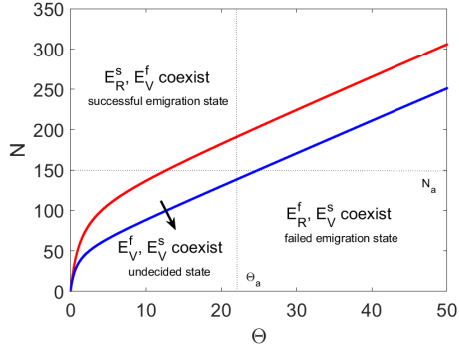}
		\caption{Bifurcation diagrams of System (3.2) with varying $N$ and $\Theta$ when $\frac{\alpha_{cs}}{\alpha_{ls}}>\max\{1,\frac{\beta_{cs}}{\beta_{ls}}\}$. Blue line is $\mathcal{N}_1(\Theta)$, red line is $\mathcal{N}_2(\Theta)$.}
		\label{Twoparametersbifurcation3}
	\end{figure}

\begin{figure}[htbp]
		\centering
		\subfloat[ ]{\label{bifurcationB:subfig1}\includegraphics[width=6cm]{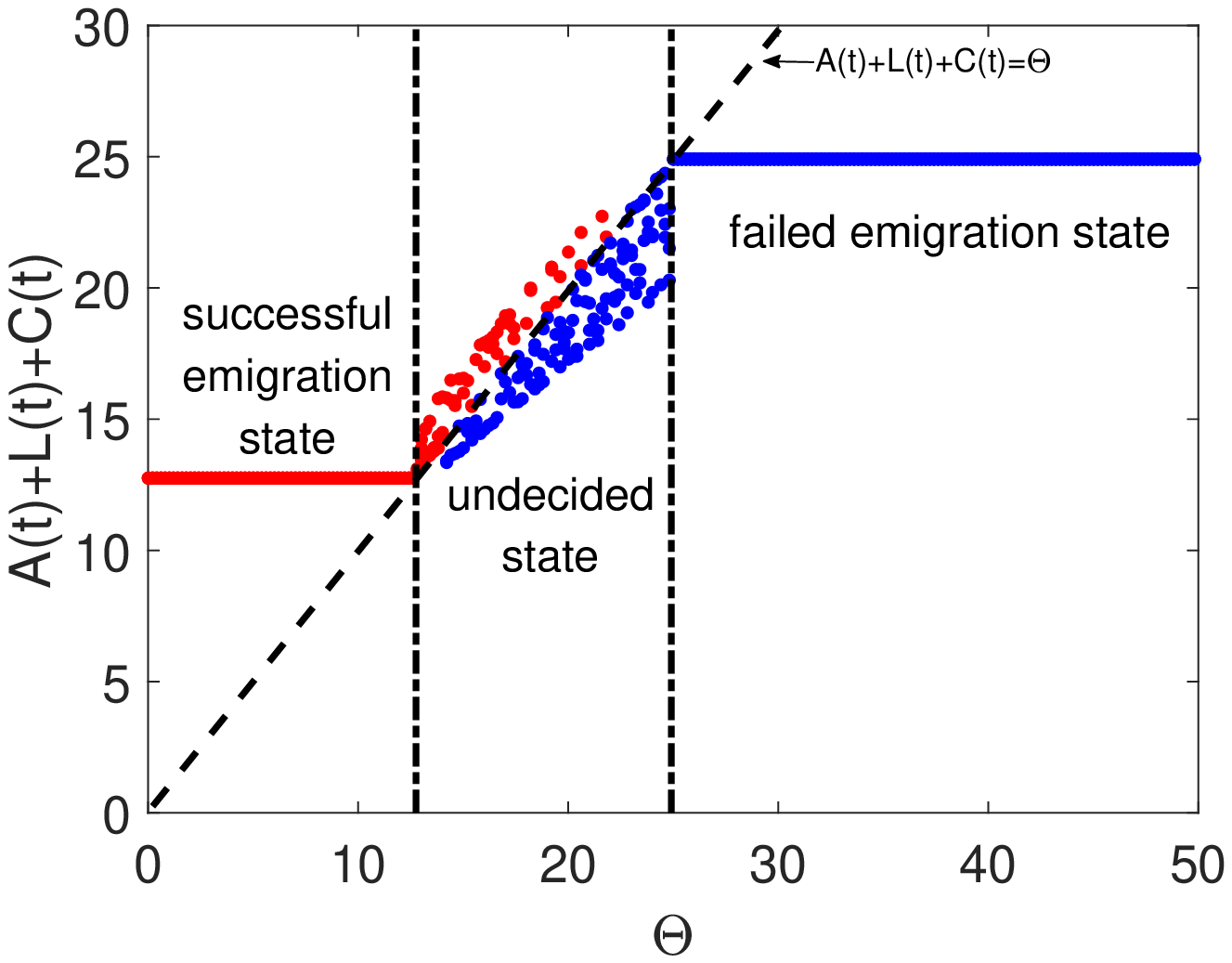}}
		\subfloat[ ]{\label{bifurcationB:subfig2}\includegraphics[width=6cm]{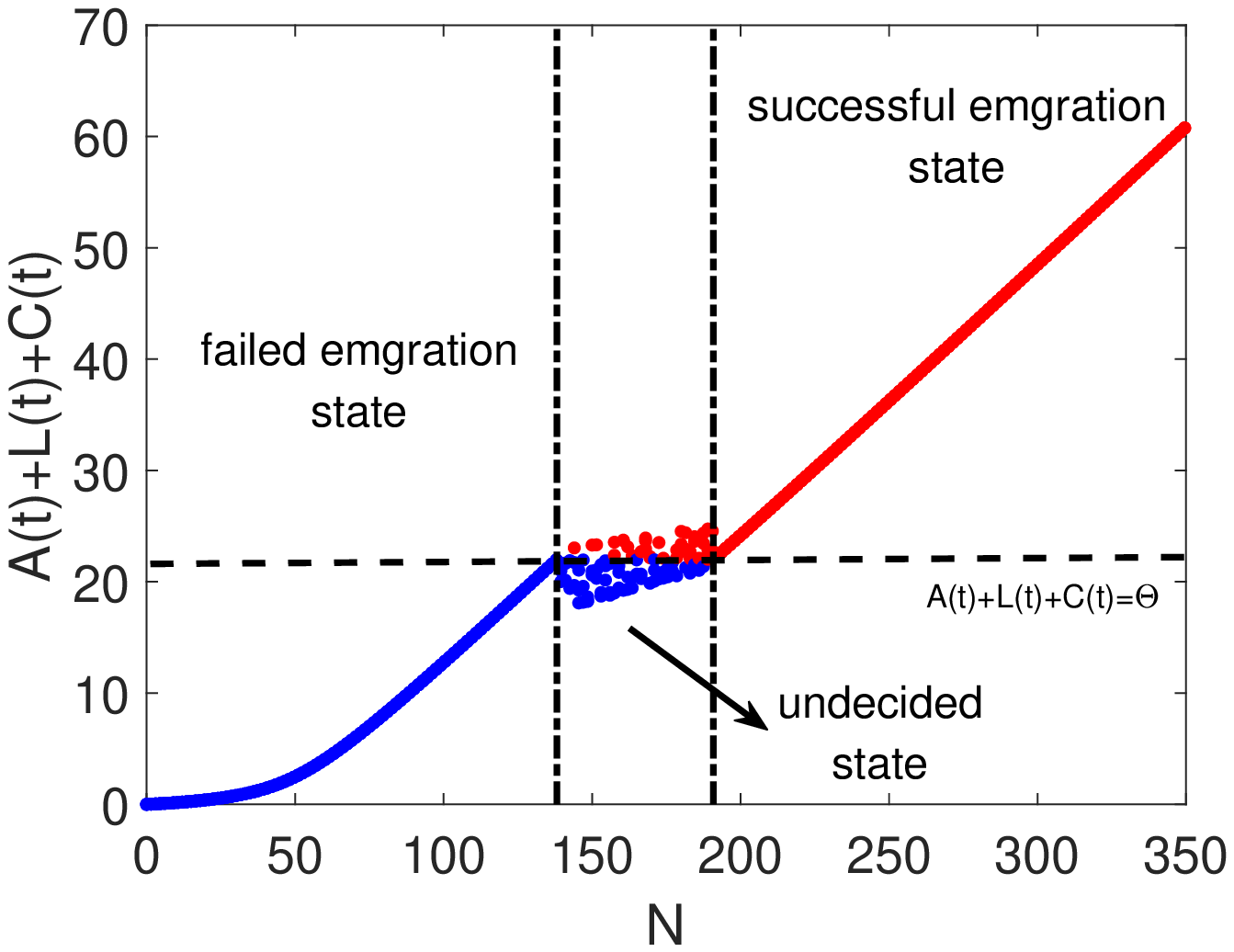}}
		\caption{Bifurcation diagrams of System (3.2). In Figure \ref{bifurcationB:subfig1}, $N=150$; In Figure \ref{bifurcationB:subfig2}, $\Theta=22$. The other parameters are $\rho=0.25$, $\alpha_{cs}=0.07$, $\alpha_{ls}=0.018$, $\beta_{ls}=0.049$, $\beta_{cs}=0.079$, $\alpha_{al}=0.007$, $\alpha_{lc}=0.15$, $\alpha_{as}=0.24$, $\alpha_{sa}=0.01$. }
		\label{bifurcationB}
	\end{figure}

	We perform bifurcation study of System (3.2) satisfying $\frac{\alpha_{cs}}{\alpha_{ls}}>\max\{1,\frac{\beta_{cs}}{\beta_{ls}}\}$. We fix $N=N_a$ in Figure \ref{Twoparametersbifurcation3} and vary $\Theta$ to obtain bifurcation diagrams as shown in Figure \ref{bifurcationB:subfig1}, and fix $\Theta=\Theta_a$ in Figure \ref{Twoparametersbifurcation3} and vary $N$ to obtain bifurcation diagrams as shown in Figure \ref{bifurcationB:subfig2}.
	
	For a colony with fixed colony size (see Figure \ref{bifurcationB:subfig1}	), when quorum threshold  is small (e.g., $\Theta$ varies from $0$ to $12$), System (3.2) stabilizes at $E_{R}^{s}$ (\textit{successful emigration state}); when quorum threshold is moderate (e.g., $\Theta$ varies from $12$ to $25$), points with two colors distribute discretely near the quorum threshold $\Theta$ i.e., System (3.2) is in undecided state; when quorum threshold is large (e.g., $\Theta$ varies from $25$ to $50$), System (3.2) stabilizes at $E_{R}^{f}$ (\textit{failed emigration state}). 
	
	For a fixed quorum threshold (see Figure \ref{bifurcationB:subfig2}), when colony size is small (e.g., $N$ varies from $0$ to $140$), System (3.2) stabilizes at $E_{R}^{f}$ (\textit{failed emigration state}); when colony size is moderate (e.g., $N$ varies from $140$ to $190$), System (3.2) is in undecided state; when colony size is large (e.g., $N$ varies from $190$ to $350$), System (3.2) stabilizes at $E_{R}^{s}$ (\textit{successful emigration state}). 
	
	\section{Bifurcation analysis when $1>\frac{\alpha_{cs}}{\alpha_{ls}}>\frac{\beta_{cs}}{\beta_{ls}}$}\label{statementD}
	
\begin{figure}[H]
		\centering
		\includegraphics[width=6.5cm]{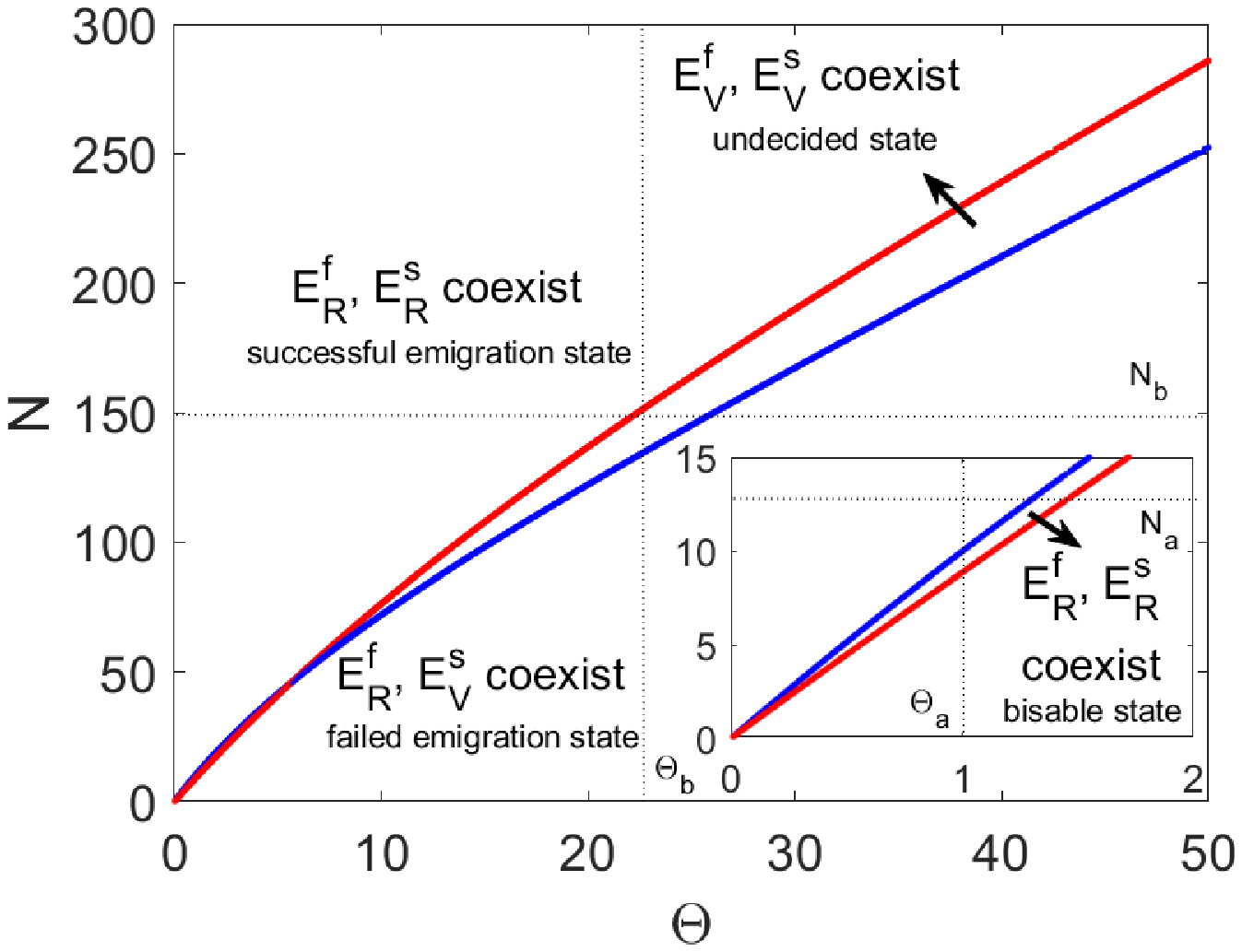}
		\caption{Bifurcation diagrams of System (3.2) with varying $N$ and $\Theta$ when $1>\frac{\alpha_{cs}}{\alpha_{ls}}>\frac{\beta_{cs}}{\beta_{ls}}$.}
		\label{Twoparametersbifurcation4}
	\end{figure}

\begin{figure}[H]
		\centering
		\subfloat[ ]{\label{bifurcationNThetaD:subfig1}\includegraphics[width=6cm]{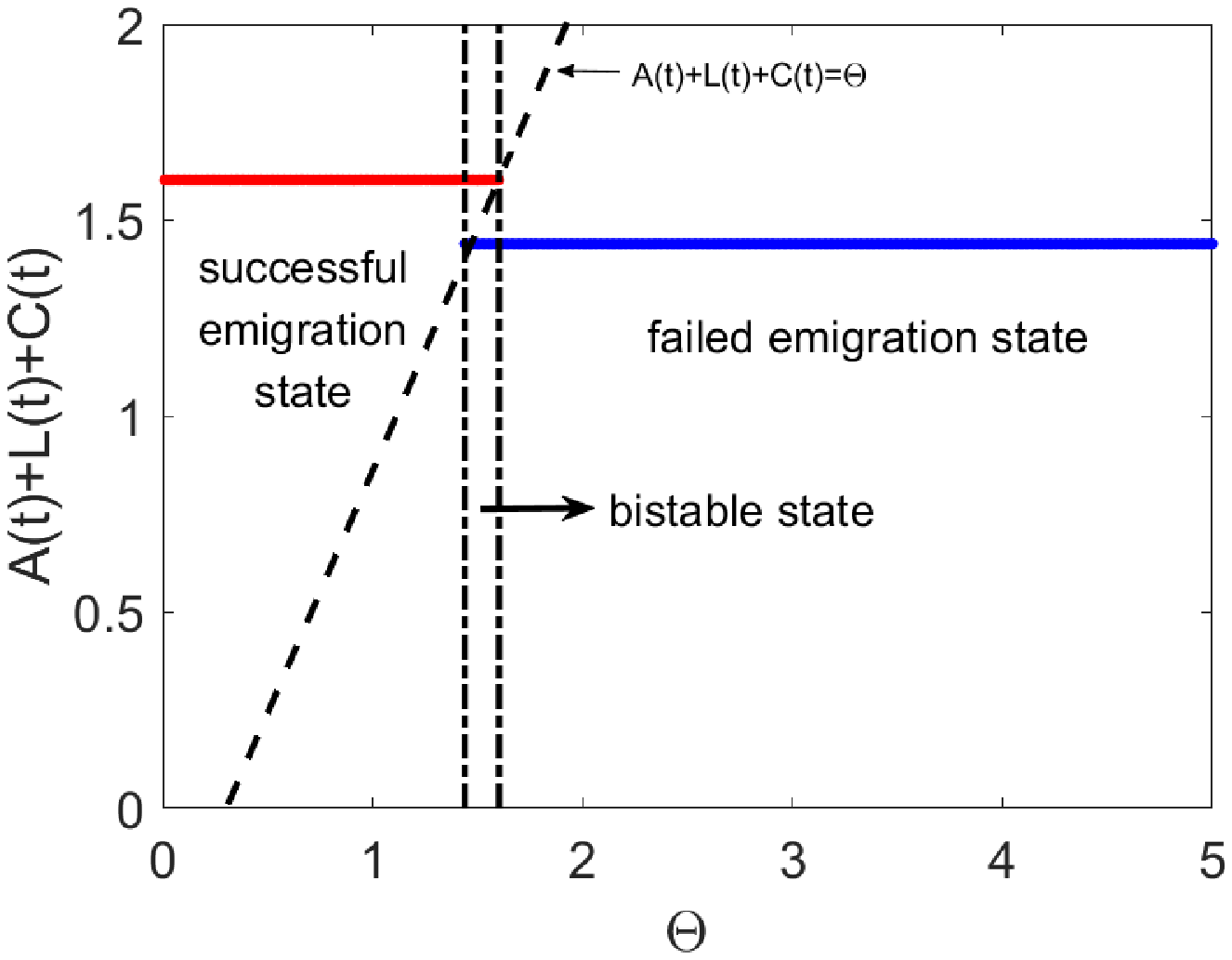}}
		\subfloat[ ]{\label{bifurcationNThetaD:subfig2}\includegraphics[width=6cm]{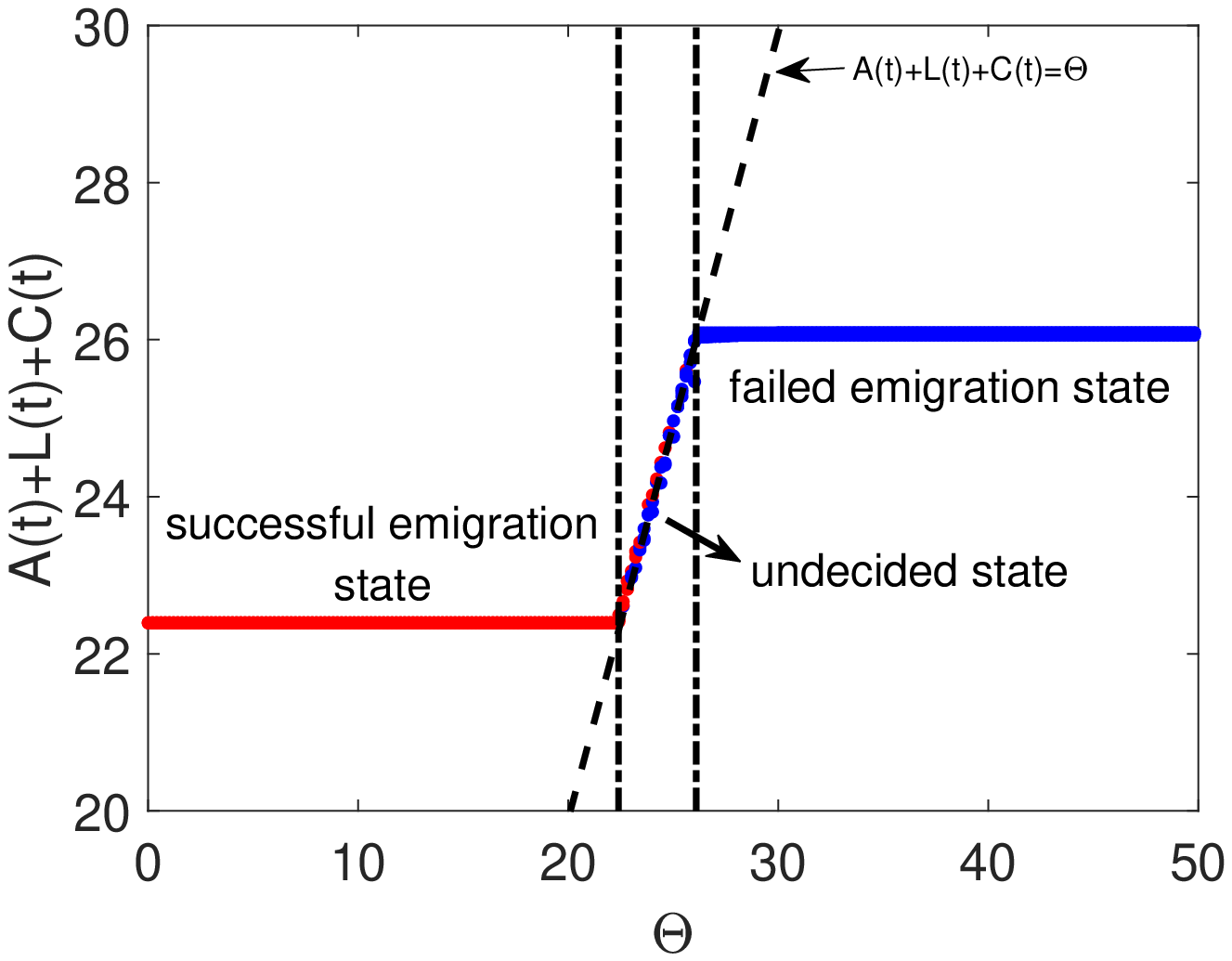}}\hspace{1mm}
		\subfloat[ ]{\label{bifurcationNThetaD:subfig3}\includegraphics[width=6cm]{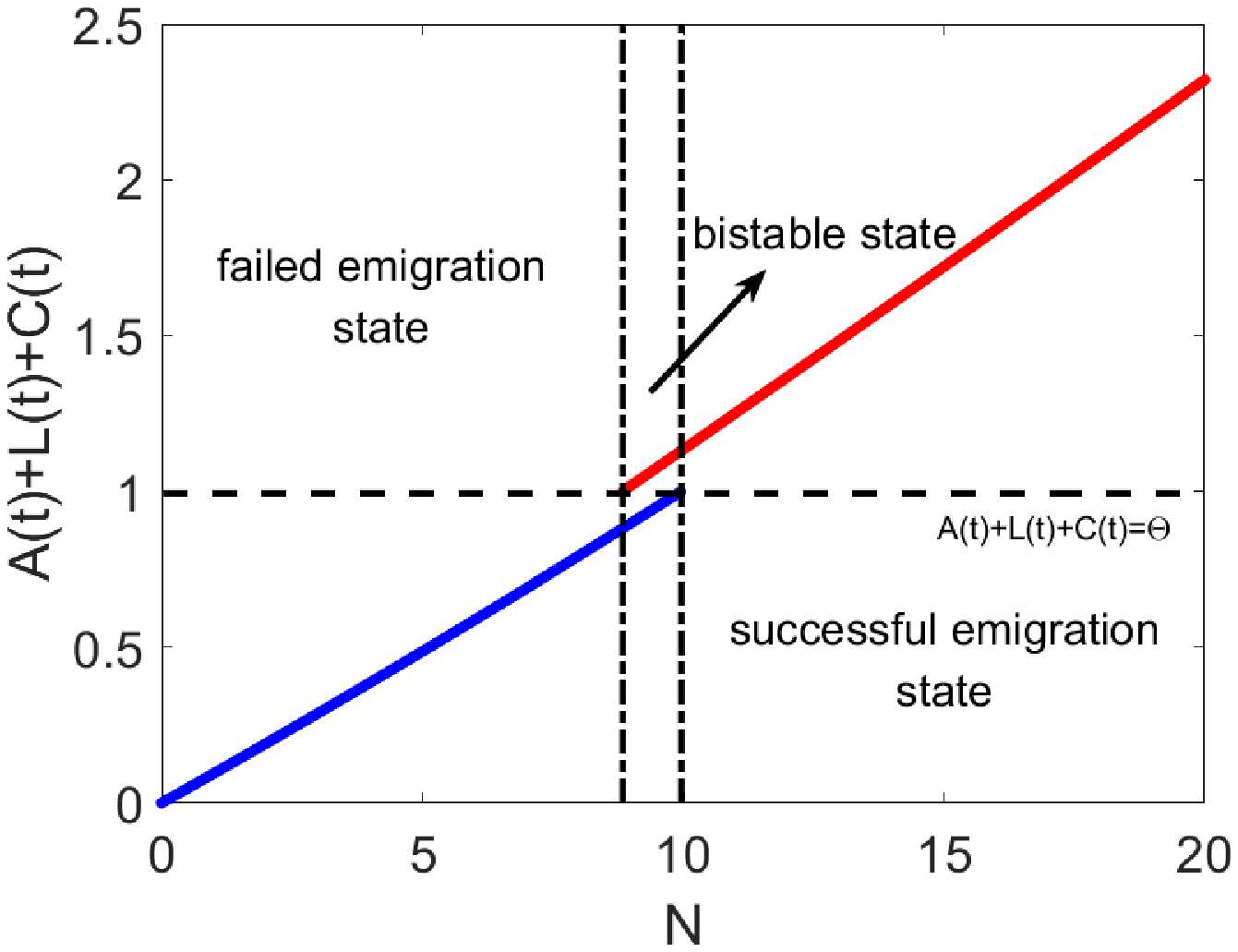}}
		\subfloat[ ]{\label{bifurcationNThetaD:subfig4}\includegraphics[width=6cm]{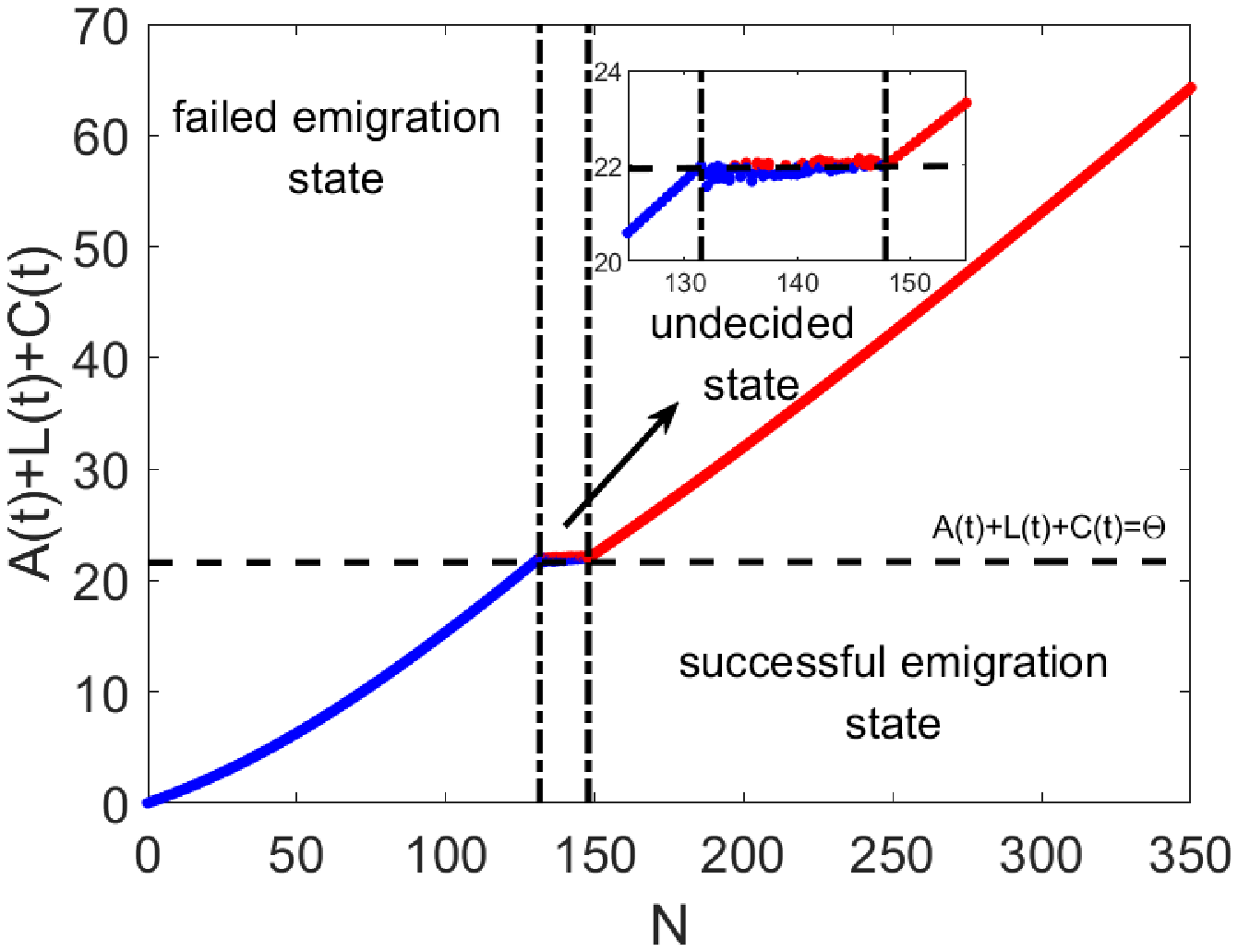}}
		\caption{Bifurcation diagrams of System (3.2) with two different levels of $N$ and two different levels of $\Theta$. In Figure \ref{bifurcationNThetaD:subfig1}, $N=14$; In Figure \ref{bifurcationNThetaD:subfig2}, $N=150$; In Figure \ref{bifurcationNThetaD:subfig3}, $\Theta=1$; In Figure \ref{bifurcationNThetaD:subfig4}, $\Theta=22$. The other parameters are $\rho=0.25$, $\alpha_{cs}=0.07$, $\alpha_{ls}=0.12$, $\beta_{ls}=0.017$, $\beta_{cs}=0.0025$, $\alpha_{al}=0.2$, $\alpha_{lc}=0.28$, $\alpha_{as}=0.24$, $\alpha_{sa}=0.1$ . }
		\label{bifurcationNThetaD}
	\end{figure}

	We perform bifurcation study of System (3.2) satisfying $1>\frac{\alpha_{cs}}{\alpha_{ls}}>\frac{\beta_{cs}}{\beta_{ls}}$. We fix two different levels of $N$ (see $N_a$ and $N_b$ in Figure \ref{Twoparametersbifurcation4}) and vary $\Theta$ to obtain bifurcation diagrams as shown in Figure \ref{bifurcationNThetaD:subfig1} and Figure \ref{bifurcationNThetaD:subfig2}, and fix two different levels of $\Theta$ (see $\Theta_a$ and $\Theta_b$ in Figure \ref{Twoparametersbifurcation4}) and vary $N$ to obtain bifurcation diagrams as shown in Figure \ref{bifurcationNThetaD:subfig3} and Figure \ref{bifurcationNThetaD:subfig4}.
	
	For a colony with small level of size (see Figure \ref{bifurcationNThetaD:subfig1}), when quorum threshold  is small (e.g., $\Theta$ varies from $0$ to $1.5$), System (3.2) stabilizes at $E_{R}^{s}$ (\textit{successful emigration state}); when quorum threshold is moderate (e.g., $\Theta$ varies from $1.5$ to $1.7$), System (3.2) is in bistability between $E_{R}^{f}$ and $E_{R}^{s}$; when quorum threshold is large (e.g., $\Theta$ varies from $1.7$ to $5$), System (3.2) stabilizes at $E_{R}^{f}$ (\textit{failed emigration state}). For a colony with large level of size (see Figure \ref{bifurcationNThetaD:subfig2}), as $\Theta$ increases, the steady-state of colony also undergoes from successful emigration state (e.g., $\Theta$ varies from $0$ to $22$) to failed emigration state (e.g., $\Theta$ varies from $26$ to $50$) but with undecided state as an intermediate (e.g., $\Theta$ varies from $22$ to $26$).

		For small quorum threshold (see Figure \ref{bifurcationNThetaD:subfig3}), when the colony size is small (e.g., $N$ varies from $0$ to $8$), System (3.2) stabilizes at $E_{R}^{f}$ (\textit{failed emigration state}); when the colony size is moderate (e.g., $N$ varies from $8$ to $10$), System (3.2) stabilizes is in bistable state between $E_{R}^{f}$ and $E_{R}^{s}$; when the colony size is large (e.g., $N$ varies from $10$ to $20$), System (3.2) stabilizes at $E_{R}^{s}$ (\textit{successful emigration state}). For large quorum threshold (see Figure \ref{bifurcationNThetaD:subfig4}), as $N$ increases, the steady-state of colony also undergoes from failed emigration state (e.g., $N$ varies from $0$ to $132$) to successful emigration state (e.g., $N$ varies from $148$ to $350$) but with undecided state as an intermediate (e.g., $N$ varies from $132$ to $148$).

	\section{Proof of Theorem 3.1}\label{proof3.1}
	\begin{proof}
		For System (2.1), we have
		\begin{displaymath}
		\frac{dS}{dt}\bigg{|}_{S=0}\geq0, \quad
		\frac{dA}{dt}\bigg{|}_{A=0}\geq0, \quad \frac{dL}{dt}\bigg{|}_{L=0}\geq0, \quad \frac{dC}{dt}\bigg{|}_{C=0}\geq0, \quad \frac{dP}{dt}\bigg{|}_{P=0}\geq0.
		\end{displaymath}
		Let $M=S+A+L+C$, then we have 
		\begin{displaymath}
		\frac{dM}{dt}=0.
		\end{displaymath}
		It then follows from $\frac{dP}{dt}=\beta_{cs}C[(1-\rho)N-P]$ that
		\begin{displaymath}
		\limsup\limits_{t\rightarrow\infty}P(t)=(1-\rho)N.
		\end{displaymath}
		Therefore, all the trajectories of System (2.1) from  $\mathbb{R}^{5}_{+}$ will enter and remain in the region \begin{displaymath}
		\Omega=\left\{(S, A, L, C, P)\in\mathbb{R}_{+}^{5}:S+A+L+C=\rho N, 0\leq P\leq (1-\rho)N \right\}.
		\end{displaymath}
		From System (2.1), there exists $\beta=\max\{\beta_{ls}, \beta_{cs}\}$, $\alpha=\min\{\alpha_{as}, \alpha_{ls}, \alpha_{cs}\}$ such that
		\begin{displaymath}
		\begin{aligned}
		\frac{dS}{dt}&=-\alpha_{sa} S-\beta_{ls}SL-\beta_{cs} SC+\alpha _{as}A+\alpha _{ls}L+\alpha _{cs}C,\\
		&\geq -\left(\alpha_{sa}+\alpha+\beta\rho N\right)S+\alpha \rho N,
		\end{aligned}
		\end{displaymath}
		which implies $\epsilon\leq\liminf\limits_{t\rightarrow\infty}S(t)\leq\limsup\limits_{t\rightarrow\infty}S(t)\leq \rho N$ with $\epsilon=\rho N\left(1-\frac{\frac{\alpha_{sa}}{\alpha}+\frac{\beta\rho N}{\alpha}}{\frac{\alpha_{sa}}{\alpha}+\frac{\beta\rho N}{\alpha}+1}\right)$. Therefore, $S$ is uniformly persistent.
		From System (2.1), we also have
		\begin{displaymath}
		\begin{aligned}
		\frac{dA}{dt}&\geq \alpha_{sa}\epsilon-\left(\alpha_{as}+\alpha_{al}\right)A,
		\end{aligned}
		\end{displaymath}
		which implies that $\liminf\limits_{t\rightarrow\infty}A(t)\geq\frac{\alpha_{sa}\epsilon}{(\alpha _{as}+\alpha_{al})}$. Letting $\epsilon_{A}=\frac{\alpha_{sa}\epsilon}{(\alpha _{as}+\alpha_{al})}$, then we have
		\begin{displaymath}
		\begin{aligned}
		\frac{dL}{dt}&\geq \alpha_{al}\epsilon_{A}-\left(\alpha_{ls}+\alpha_{lc}\right)L.
		\end{aligned}
		\end{displaymath}
		It then follows that $\liminf\limits_{t\rightarrow\infty}L(t)\geq\frac{\alpha_{al}\epsilon_{A}}{(\alpha _{ls}+\alpha_{lc})}$.
		Thus, the persistence of $S$ leads to the persistence of $A$ and $L$.
	\end{proof}
	
	\section{Proof of Theorem 3.2}
	We first denote function 
	\begin{displaymath}
	g_1(L)=-\beta_{ls}(1+\frac{\alpha_{ls}}{\alpha_{al}})L^2+\left[\beta_{ls}\rho N-\alpha_{sa}\left(1+\frac{\alpha_{ls}}{\alpha_{al}}\right)-(\alpha_{as}+\alpha_{al})\frac{\alpha_{ls}}{\alpha_{al}}\right]L+\alpha_{sa}\rho N
	\end{displaymath}
	with one positive root
	\begin{displaymath}
	L^{*}=\frac{\left(\rho N \xi_1-1-\eta_1\right)+\sqrt{\left(\rho N \xi_1-1-\eta_1\right)^2+4\xi_1\rho N}}{2\xi_1\left(1+\frac{\alpha_{ls}}{\alpha_{al}}\right)}.
	\end{displaymath}
	Then, Subsystem (3.3) has unique equilibrium $E^{f}(A^{f}, L^{f}, 0)$ where $L^{f}=L^{*}, A^{f}=\frac{\alpha_{ls}}{\alpha_{al}}L^{f}$.
	
	The Jacobin matrix of subsystem (3.3) at $E^{f}$ is presented as follows
	\begin{displaymath}
	J(E^{f})=\left[
	\begin{array}{cc}
	-(\alpha_{sa}+\beta_{ls} L^f)-(\alpha_{as}+\alpha_{al})& \beta_{ls}(\rho N-A^f-L^f)-(\alpha_{sa}+\beta_{ls} L^f) \\
	\alpha_{al} & -\alpha_{ls}
	\end{array}
	\right]. 
	\end{displaymath}
	After extensive algebraic calculations, the corresponding characteristic equation at $E^{f}$ is
	\begin{equation}\label{characteristic1}
	\lambda^2+a_{1}\lambda+a_{0}=0,
	\end{equation}
	where
	\begin{displaymath}
	\begin{aligned}
	a_{1}&=(\alpha_{as}+\alpha_{al})+(\alpha_{sa}+\beta_{ls} L^{f})+\alpha_{ls},\\
	a_{0}&=2\beta_{ls}(\alpha_{ls}+\alpha_{al})L^{f}+\alpha_{sa}(\alpha_{ls}+\alpha_{al})+\alpha_{ls}(\alpha_{as}+\alpha_{al})-\alpha_{al}\beta_{ls}\rho N,\\
	&=-\alpha_{al}g_1'(L^{f}).
	\end{aligned}
	\end{displaymath}
	Let $\lambda_{1}(E^{f})$ and $\lambda_{2}(E^{f})$ be the roots of \eqref{characteristic1} with $\Re\lambda_{1}(E^{f})\leq\Re\lambda_{2}(E^{f})$. It then follows that
	\begin{displaymath}
	\begin{aligned}
	\lambda_{1}(E^{f})+\lambda_{2}(E^{f})&=-a_1<0,\\
	\lambda_{1}(E^{f})\lambda_{2}(E^{f})&=a_0>0.
	\end{aligned}
	\end{displaymath}
	Thus, the positive equilibrium $E^{f}$ is locally asymptotically stable.  Next, we prove that $E^{f}$ is globally asymptotically stable.
	
	From subsystem (3.3), it follows that 
	\begin{equation}\label{eq:1}
	\lim\limits_{t\rightarrow+\infty}C(t)=0.
	\end{equation}
	Substituting \eqref{eq:1} into the subsystem (3.3), we obtain the limit system
	\begin{equation}\label{limitsystem}
	\begin{aligned}
	\frac{dA}{dt}&=\left(\alpha_{sa}+\beta_{ls}L\right)\left(\rho N-A-L\right)-\alpha _{as}A-\alpha_{al}  A,\\
	\frac{dL}{dt}&=\alpha_{al}  A-\alpha _{ls}L,\\
	\end{aligned}
	\end{equation}
	From the Jacobian of subsystem (3.3) at equilibrium $E^{f}(A^{f}, L^{f}, 0)$, we can easily get that $(A^{f}, L^{f})$ is locally stable.	Let $h_1(A, L)=(\alpha_{sa}+\beta_{ls}L)(\rho N-A-L)-(\alpha_{al}+\alpha_{as})A$ and $h_2(A, L)=\alpha_{al} A-\alpha_{ls} L.$ It then follows that
	\begin{displaymath}
	\frac{\partial{h_1(A,L)}}{\partial{A}}+\frac{\partial{h_2(A,L) }}{\partial{L}}=-(\alpha_{al}+\alpha_{as}+\alpha_{ls}+\alpha_{sa})-\beta_{ls}L<0 \quad \text{for all $L>0$.}
	\end{displaymath}
	By Poincare-Bendixson Theorem, there is no limit cycle in the system \eqref{limitsystem}. Since the omega limit set of the system \eqref{limitsystem} is either a fixed point or a limit cycle, the interior equilibrium $(A^{f}, L^{f})$ is globally stable. Since system \eqref{limitsystem} is the limit system of (3.3), it follows that $E^{f}(A^{f}, L^{f}, 0)$ is globally asymptotically stable for system (3.3). 
	
	By using the same argument, we denote function
	\begin{displaymath}
	\begin{aligned}
	g_2(L)=&-\left(\beta_{ls}\frac{\alpha_{cs}}{\alpha_{lc}}+\beta_{cs}\right)\left[1+\frac{\alpha_{cs}}{\alpha_{lc}}+\frac{\alpha_{cs}(\alpha_{lc}+\alpha_{ls})}{\alpha_{al}\alpha_{lc}}\right]C^2\\
	&-\left[\alpha_{sa}\left(1+\frac{\alpha_{cs}}{\alpha_{lc}}+\frac{\alpha_{cs}(\alpha_{lc}+\alpha_{ls})}{\alpha_{al}\alpha_{lc}}\right)+\frac{\alpha_{cs}(\alpha_{lc}+\alpha_{ls})(\alpha_{as}+\alpha_{al})}{\alpha_{al}\alpha_{lc}}\right]C\\
	&+\rho N\left(\beta_{ls}\frac{\alpha_{cs}}{\alpha_{lc}}+\beta_{cs}\right)C+\alpha_{sa}\rho N,
	\end{aligned}
	\end{displaymath}
	and \begin{displaymath}
	C^{*}=\frac{\left(\rho N \xi_2-1-\eta_2\right)+\sqrt{\left(\rho N \xi_2-1-\eta_2\right)^2+4\xi_2\rho N}}{2\xi_2\left[1+\frac{\alpha_{cs}}{\alpha_{lc}}+\frac{\alpha_{cs}\left(\alpha_{lc}+\alpha_{ls}\right)}{\alpha_{al}\alpha_{lc}}\right]}
	\end{displaymath} is the unique positive root of $g_2(L)=0$. Then, it is easy to get that Subsystem (3.4) has unique positive equilibrium $E^{s}(A^{s}, L^{s},C^{s})$, where 
	\begin{displaymath}
	C^{s}=C^{*}, \quad  L^{s}=\frac{\alpha_{cs}}{\alpha_{lc}}C^{s}, \quad  A^{s}=\frac{\alpha_{cs}(\alpha_{lc}+\alpha_{ls})}{\alpha_{al}\alpha_{lc}}C^{s}.
	\end{displaymath}

	The Jacobin matrix of subsystem (3.4) at $E^{s}$ is presented as follows
	\begin{displaymath}
	J(E^s)=\left[
	\begin{array}{ccc}
	A_{11}& A_{12} & A_{13}\\
	\alpha_{al} & -(\alpha_{lc}+\alpha_{ls})  & 0 \\
	0 & \alpha_{lc} & -\alpha_{cs} \\
	\end{array}
	\right], 
	\end{displaymath}
	where $$
	\begin{aligned}
	A_{11}=&-(\alpha_{sa}+\alpha_{as}+\alpha_{al})-\beta_{ls} L^{s}-\beta_{cs} C^{s},\\
	A_{12}=&(\beta_{ls}\rho N-\alpha_{sa})-\beta_{ls}A^{s}-2\beta_{ls}L^{s}-(\beta_{ls}+\beta_{cs})C^{s},\\
	A_{13}=&(\beta_{cs}\rho N-\alpha_{sa})-\beta_{cs}A^{s}-(\beta_{ls}+\beta_{cs})L^{s}-2\beta_{cs}C^{s}.
	\end{aligned}
	$$
	After extensive algebraic calculations, the corresponding characteristic equation at $E^{s}$ is
	\begin{equation}\label{characteristic2}
	\lambda^3+b_{2}\lambda^2+b_{1}\lambda+b_{0}=0,
	\end{equation}
	where
	\begin{displaymath}
	\begin{aligned}
	b_{2}=&(\alpha_{cs}+\alpha_{ls}+\alpha_{as}+\alpha_{al}+\alpha_{lc}+\alpha_{sa})+\beta_{ls}L^{s}+\beta_{cs}C^{s},\\
	b_{1}=&\alpha_{cs}(\alpha_{ls}+\alpha_{as}+\alpha_{al}+\alpha_{lc}+\alpha_{sa}+\beta_{ls}L^{s}+\beta_{cs}C^{s})\\&+(\alpha_{al}+\alpha_{lc}+\alpha_{ls})(\alpha_{sa}+\beta_{ls}L^{s}+\beta_{cs}C^{s})\\
	&+(\alpha_{sa}+\beta_{cs}C^{s})(\rho N-A^{s}-L^{s}-C^{s})\frac{\alpha_{al}}{L^{s}},\\
	b_{0}=&-\alpha_{al}\alpha_{lc} g_2'(C^{s}).
	\end{aligned}
	\end{displaymath}
	Let $\lambda_{1}(E^{s})$, $\lambda_{2}(E^{s})$ and $\lambda_{3}(E^{s})$ be the roots of \eqref{characteristic2}. It then follows that
	\begin{equation}\label{characteristic3}
	\begin{aligned}
	\lambda_{1}(E^{s})+\lambda_{2}(E^{s})+\lambda_{3}(E^{s})&=-b_2<0,\\
	\lambda_{1}(E^{s})\lambda_{2}(E^{s})\lambda_{3}(E^{s})&=-b_0<0.
	\end{aligned}
	\end{equation}
	Assume that $\Re\lambda_{1}(E^{s})\leq\Re\lambda_{2}(E^{s})\leq\Re\lambda_{3}(E^{s})$. Then, \eqref{characteristic3} implies that $\Re\lambda_{1}(E^{s})\leq\Re\lambda_{2}(E^{s})\leq\Re\lambda_{3}(E^{s})<0$ or $\Re\lambda_{1}(E^{s})<0<\Re\lambda_{2}(E^{s})\leq\Re\lambda_{3}(E^{s})$. Simple calculations give that $b_{1}b_{2}-b_{0}>0$ holds for all parameters. Thus, $\Re\lambda_{1}(E^{s})\leq\Re\lambda_{2}(E^{s})\leq\Re\lambda_{3}(E^{s})<0$. Based on the Routh-Hurwitz criterion, the positive equilibrium $E^{s}$ is locally asymptotically stable.
	
	Using the change of variables $W=A+L+C$, $V=L+C$ and $U=C$, the subsystem (3.4) can be written as
	\begin{equation}\label{transition}
	\begin{aligned}
	\frac{dW}{dt}&=\left[\alpha_{sa}+\beta_{ls}(V-U)+\beta_{cs} U\right]\left(\rho N-W\right)-\alpha _{as}(W-V)-\alpha_{ls}(V-U)-\alpha_{cs}U,\\
	\frac{dV}{dt}&=\alpha_{al}  (W-V)-\alpha _{ls}(V-U)-\alpha_{cs}U,\\
	\frac{dU}{dt}&=\alpha_{lc} (V-U)-\alpha _{cs}U.
	\end{aligned}
	\end{equation}
	It follows that the subsystem \eqref{transition} has a unique positive equilibrium $(W^{s}, V^{s}, U^{s})$ which is always locally asymptotically stable. The Jacobian has the form 
	\begin{displaymath}
	J=\left[
	\begin{array}{ccc}
	-\left[\alpha_{as}+\alpha_{sa}+\beta_{ls}V+(\beta_{cs}-\beta_{ls}) U\right]& A_{12} & A_{13}\\
	\alpha_{al} & -(\alpha_{lc}+\alpha_{ls})  & A_{23} \\
	0 & \alpha_{lc} & -(\alpha_{cs}+\alpha_{lc}) \\
	\end{array}
	\right],
	\end{displaymath}
	where $$
	\begin{aligned}
		A_{12}=&\beta_{ls}(\rho N-W)+\alpha_{as}-\alpha_{ls},\\ A_{13}=&(\beta_{cs}-\beta_{ls})(\rho N-W)+\alpha_{ls}-\alpha_{cs}, \\ A_{23}=&\alpha_{ls}-\alpha_{cs}.
	\end{aligned}
	$$ 
	
	Note that the off-diagonal entries of the Jocabian matrix are nonnegative in Int$\Omega$ if $\alpha_{as}>\alpha_{ls}>\alpha_{cs}$ and $\beta_{cs}>\beta_{ls}$. Thus, under this condition, the system \eqref{transition} is cooperative in Int$\Omega$. Furthermore, the Jocabian matrix is irreducible in Int$\Omega$. It then follows that the flow generated by \eqref{transition} is strongly monotone in Int$\Omega$ if $\alpha_{as}>\alpha_{ls}>\alpha_{cs}$ and $\beta_{cs}>\beta_{ls}$. Since the system \eqref{transition} is dissipative and the system \eqref{transition} has a unique positive equilibrium $(W^{s}, V^{s}, U^{s})$ which is always locally asymptotically stable, the equilibrium $(W^{s}, V^{s}, U^{s})$ is globally asymptotically stable. Thus, the positive equilibrium $E^{s}$ is globally asymptotically stable if $\alpha_{as}>\alpha_{ls}>\alpha_{cs}$ and $\beta_{cs}>\beta_{ls}$.
	This completes the proof.
	
	\section{Proof of Theorem 5.1} 
	Rewritten
	\begin{displaymath}
	\begin{aligned}
\mathcal{N}_{1}(\Theta)=&\frac{\xi_1\Theta^2+\Theta(1+\eta_1)}{\rho\left(1+\xi_1\Theta\right)}, \\ \mathcal{N}_{2}(\Theta)=&\frac{\xi_2\Theta^2+\Theta(1+\eta_2)}{\rho\left(1+\xi_2\Theta\right)}.
	\end{aligned}
	\end{displaymath}
	Simple calculations yield that
	\begin{displaymath}
	\begin{aligned}
\frac{d{\mathcal{N}_1(\Theta)}}{d{\Theta}}=&\frac{2\xi_1\Theta+\eta_1+1+\xi_1^2\Theta^2}{\rho(1+\xi_1\Theta)^2}>0,\\
\frac{d{\mathcal{N}_2(\Theta)}}{d{\Theta}}=&\frac{2\xi_2\Theta+1+\eta_2+\xi_2^2\Theta^2}{\rho(1+\xi_2\Theta)^2}>0,
	\end{aligned}
	\end{displaymath} and 
	\begin{displaymath}
	\begin{aligned}
	\mathcal{N}_{1}(\Theta)-\mathcal{N}_{2}(\Theta)
	=&\frac{\Theta(\xi_1\eta_2-\xi_2\eta_1)}{\rho(1+\xi_1\Theta)(1+\xi_2\Theta)}\left(\frac{\eta_1-\eta_2}{\xi_1\eta_2-\xi_2\eta_1}-\Theta\right).
	\end{aligned}
	\end{displaymath} It then follows that both $\mathcal{N}_1(\Theta)$ and $\mathcal{N}_2(\Theta)$ are monotonic increasing function on $\Theta\in(0, +\infty)$, and
	\begin{enumerate}
		\item if $\xi_1\eta_2-\xi_2\eta_1>0$ and $\eta_1-\eta_2<0$, i.e., $\frac{\alpha_{cs}}{\alpha_{ls}}<1$ and $\frac{\alpha_{cs}}{\alpha_{ls}}<\frac{\beta_{cs}}{\beta_{ls}}$, then $\mathcal{N}_{1}(\Theta)>\mathcal{N}_{2}(\Theta)$ for all $\Theta>0$;
		\item if $\xi_1\eta_2-\xi_2\eta_1<0$ and $\eta_1-\eta_2>0$, i.e., $\frac{\alpha_{cs}}{\alpha_{ls}}>1$ and $\frac{\alpha_{cs}}{\alpha_{ls}}>\frac{\beta_{cs}}{\beta_{ls}}$, then $\mathcal{N}_{1}(\Theta)<\mathcal{N}_{2}(\Theta)$ for all $\Theta>0$;
		\item if $\xi_1\eta_2-\xi_2\eta_1>0$ and $\eta_1-\eta_2>0$, i.e., $\frac{\alpha_{cs}}{\alpha_{ls}}>1$ and $\frac{\alpha_{cs}}{\alpha_{ls}}<\frac{\beta_{cs}}{\beta_{ls}}$, then $\mathcal{N}_{1}(\Theta)>\mathcal{N}_{2}(\Theta)$ for all $0<\Theta<\frac{\alpha_{sa}}{\beta_{ls}}\frac{1-\frac{\alpha_{cs}}{\alpha_{ls}}}{\frac{\alpha_{cs}}{\alpha_{ls}}-\frac{\beta_{cs}}{\beta_{ls}}}$, $\mathcal{N}_{1}(\Theta)<\mathcal{N}_{2}(\Theta)$ for all $\Theta>\frac{\alpha_{sa}}{\beta_{ls}}\frac{1-\frac{\alpha_{cs}}{\alpha_{ls}}}{\frac{\alpha_{cs}}{\alpha_{ls}}-\frac{\beta_{cs}}{\beta_{ls}}}$, and $\mathcal{N}_{1}(\Theta)=\mathcal{N}_{2}(\Theta)$ at $\Theta=\frac{\alpha_{sa}}{\beta_{ls}}\frac{1-\frac{\alpha_{cs}}{\alpha_{ls}}}{\frac{\alpha_{cs}}{\alpha_{ls}}-\frac{\beta_{cs}}{\beta_{ls}}}$;
		\item if $\xi_1\eta_2-\xi_2\eta_1<0$ and $\eta_1-\eta_2<0$, i.e., $\frac{\alpha_{cs}}{\alpha_{ls}}<1$ and $\frac{\alpha_{cs}}{\alpha_{ls}}>\frac{\beta_{cs}}{\beta_{ls}}$, then $\mathcal{N}_{1}(\Theta)<\mathcal{N}_{2}(\Theta)$ for all $0<\Theta<\frac{\alpha_{sa}}{\beta_{ls}}\frac{1-\frac{\alpha_{cs}}{\alpha_{ls}}}{\frac{\alpha_{cs}}{\alpha_{ls}}-\frac{\beta_{cs}}{\beta_{ls}}}$, $\mathcal{N}_{1}(\Theta)>\mathcal{N}_{2}(\Theta)$ for all $\Theta>\frac{\alpha_{sa}}{\beta_{ls}}\frac{1-\frac{\alpha_{cs}}{\alpha_{ls}}}{\frac{\alpha_{cs}}{\alpha_{ls}}-\frac{\beta_{cs}}{\beta_{ls}}}$ and $\mathcal{N}_{1}(\Theta)=\mathcal{N}_{2}(\Theta)$ at $\Theta=\frac{\alpha_{sa}}{\beta_{ls}}\frac{1-\frac{\alpha_{cs}}{\alpha_{ls}}}{\frac{\alpha_{cs}}{\alpha_{ls}}-\frac{\beta_{cs}}{\beta_{ls}}}$.
	\end{enumerate} This completes proof.

	\newpage
	\section{Additional Figures}
	Here, we provide the figures that can be referenced from the main text.
	\begin{figure}[htbp]
		\centering
		\includegraphics[width=14cm]{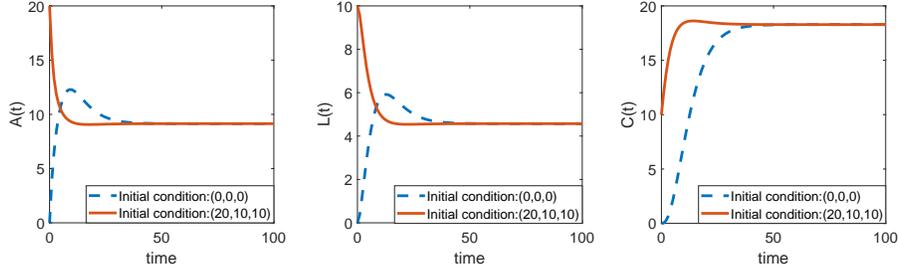}
		\caption{Time series plots of population $A$, $L$ and $C$ of subsystem (3.4). The paraeters are $N=200$, $\rho=0.25$, $\alpha_{cs}=0.07$, $\alpha_{ls}=0.12$, $\beta_{ls}=0.033$, $\beta_{cs}=0.079$, $\alpha_{al}=0.032$, $\alpha_{lc}=0.15$, $\alpha_{as}=0.24$, $\alpha_{sa}=0.01$.}
		\label{GlobalStable}
	\end{figure}
	
	\begin{figure}[htbp]
		\centering
		\subfloat[Crossing]{\label{Set:subfig1}\includegraphics[width=5cm]{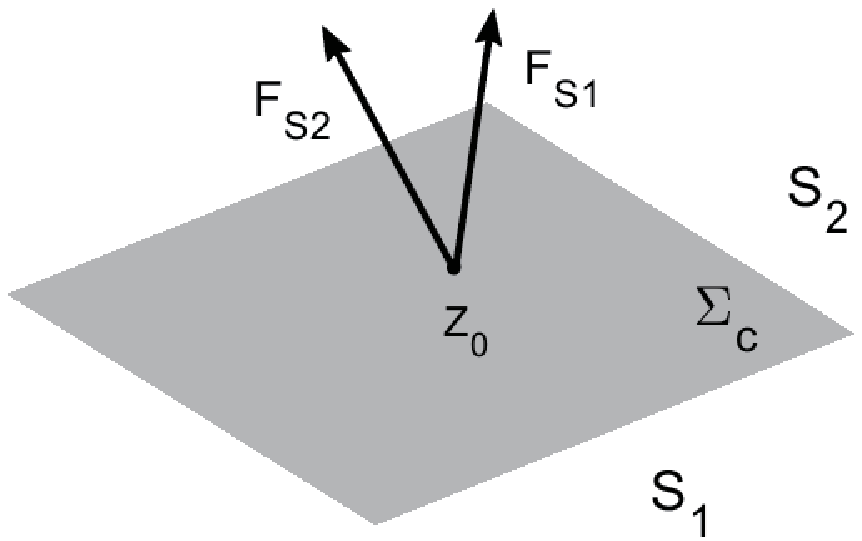}}\hspace{10mm}
		\subfloat[Sliding]{\label{Set:subfig2}\includegraphics[width=5cm]{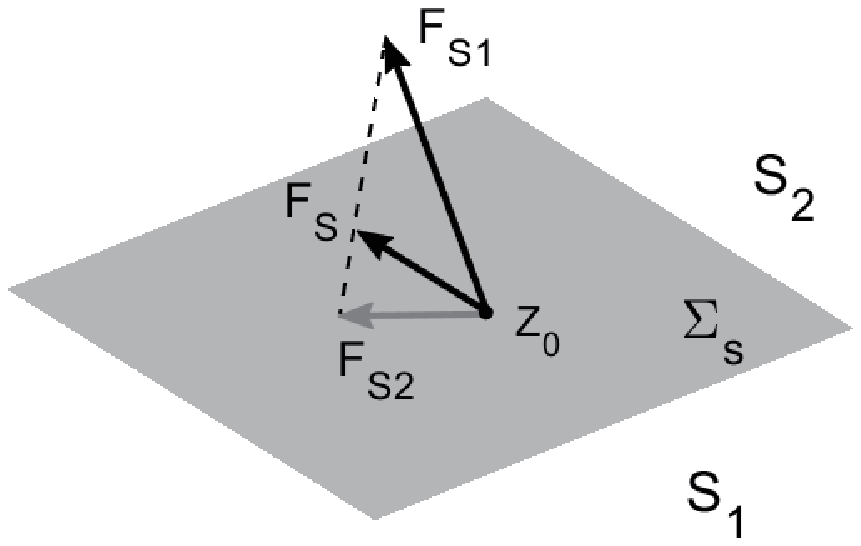}}
		\caption{Two possible vector fields on manifold $\Sigma$: (a) crossing (b) sliding along the sliding vector $F_{S}$.}
	\end{figure}
	
	\begin{figure}[htbp]
		\centering
		\subfloat[ ]{\label{ScenarioA:subfig1}\includegraphics[width=5cm]{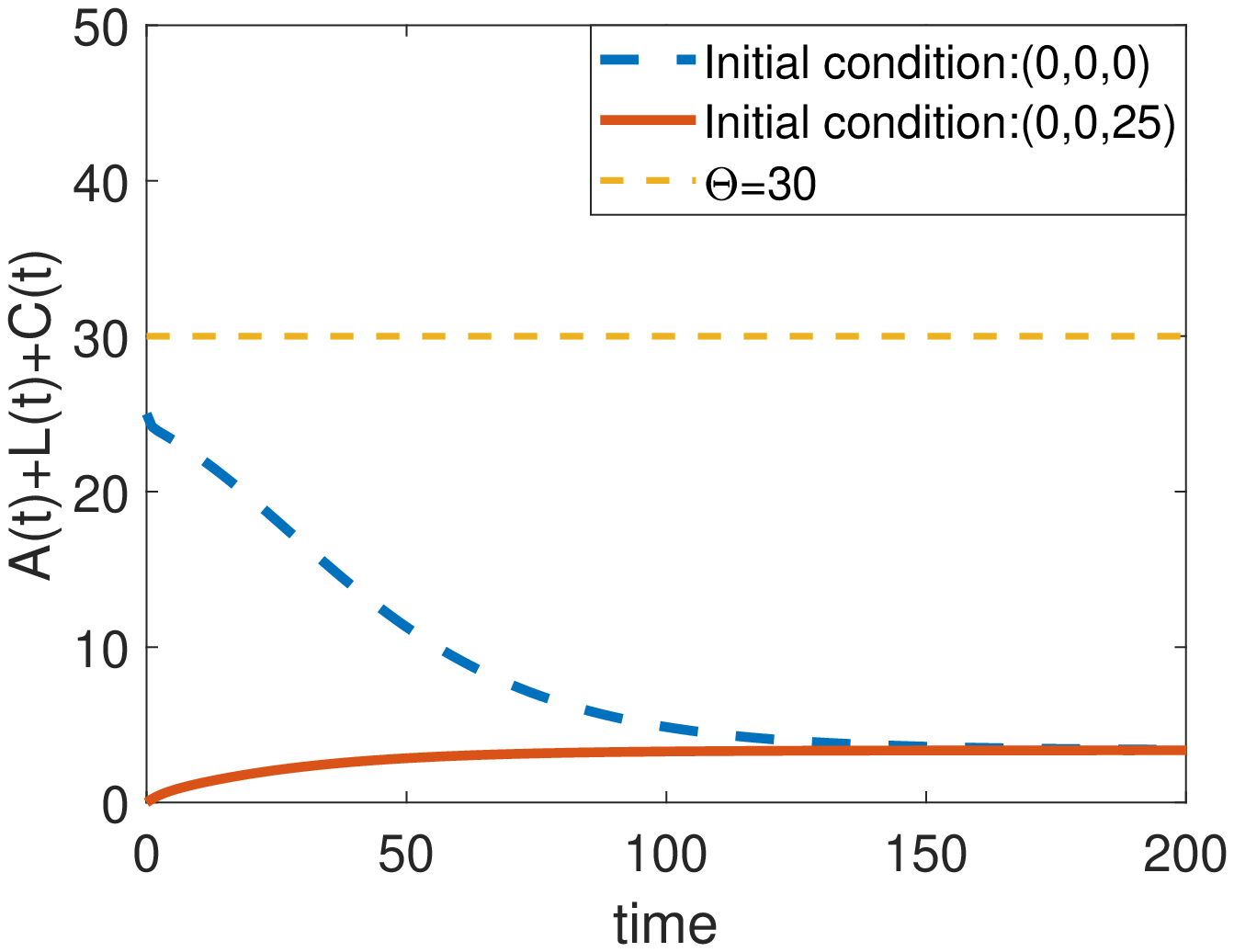}}\hspace{3mm}
		\subfloat[ ]{\label{ScenarioA:subfig2}\includegraphics[width=5.3cm]{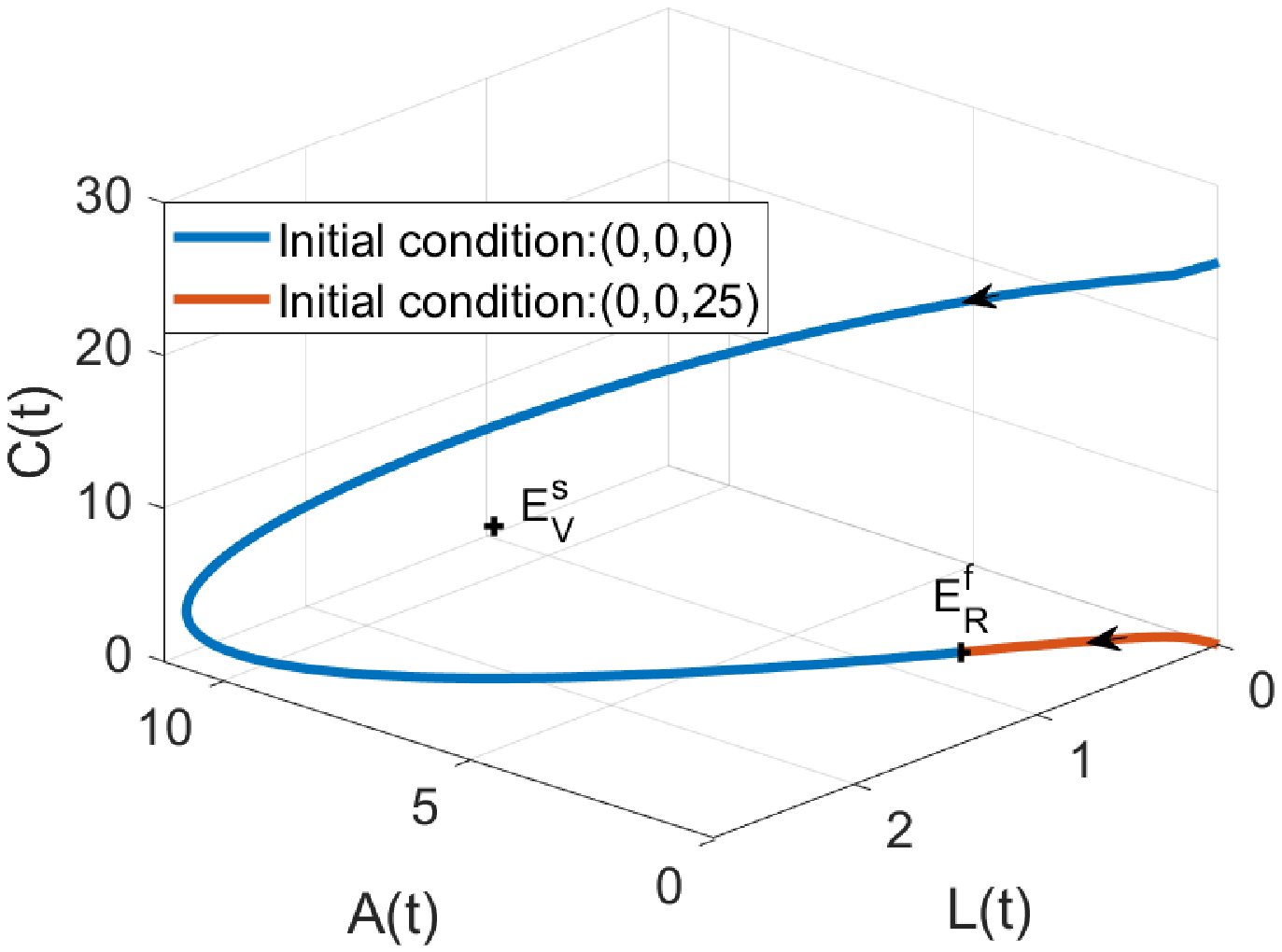}}
		\caption{Time series and phase plots of System (3.2) shows that equilibrium $E_{R}^{f}$ is the unique attractor when the parameters are $N=100$, $\Theta=30$, $\rho=0.25$, $\alpha_{cs}=0.07$, $\alpha_{ls}=0.12$, $\beta_{ls}=0.033$, $\beta_{cs}=0.079$, $\alpha_{al}=0.032$, $\alpha_{lc}=0.15$, $\alpha_{as}=0.24$, $\alpha_{sa}=0.01$.}
		\label{ScenarioA}
	\end{figure}
	
	\begin{figure}[htbp]
		\centering
		\subfloat[ ]{\label{ScenarioB:subfig1}\includegraphics[width=5cm]{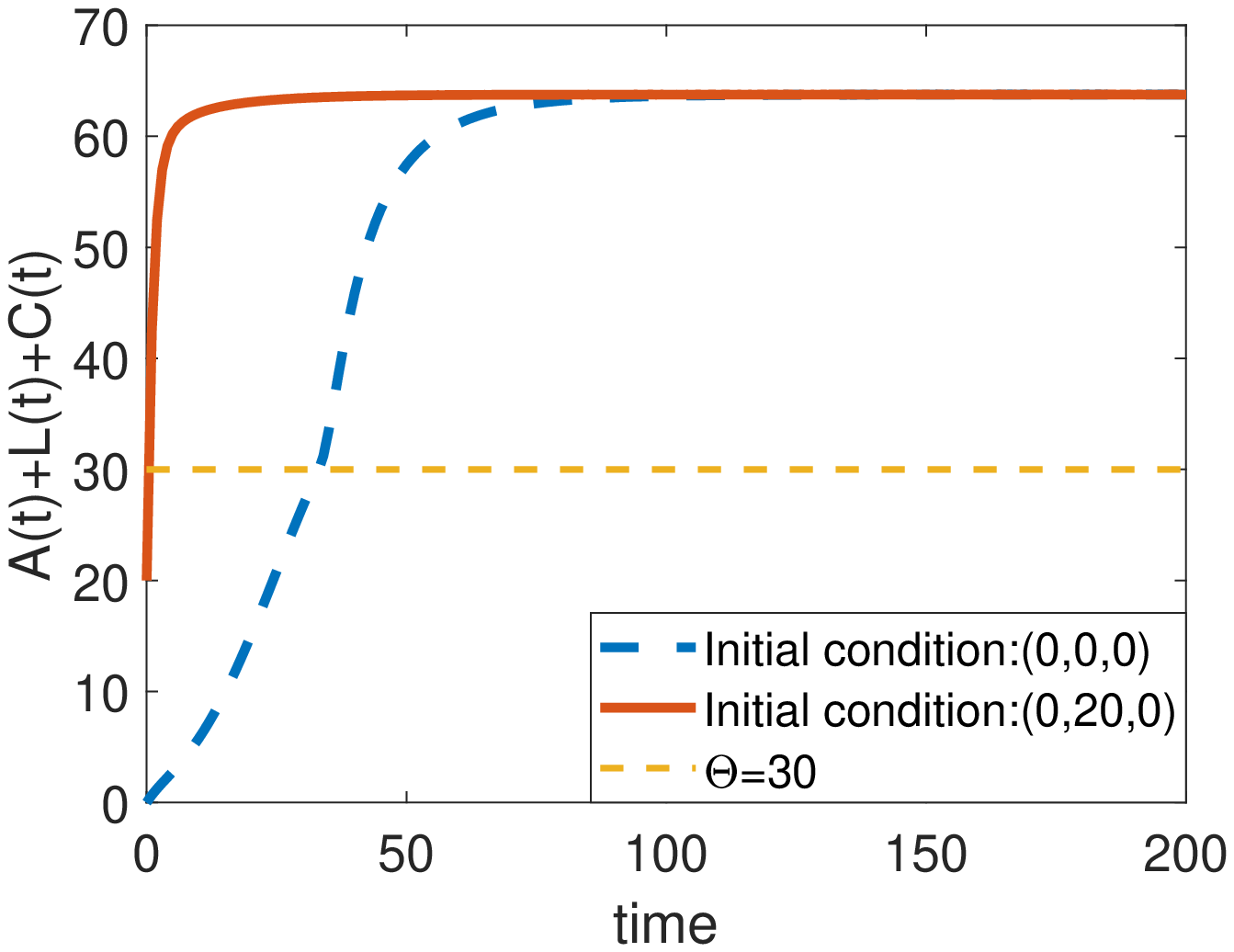}}\hspace{3mm}
		\subfloat[ ]{\label{ScenarioB:subfig2}\includegraphics[width=5.3cm]{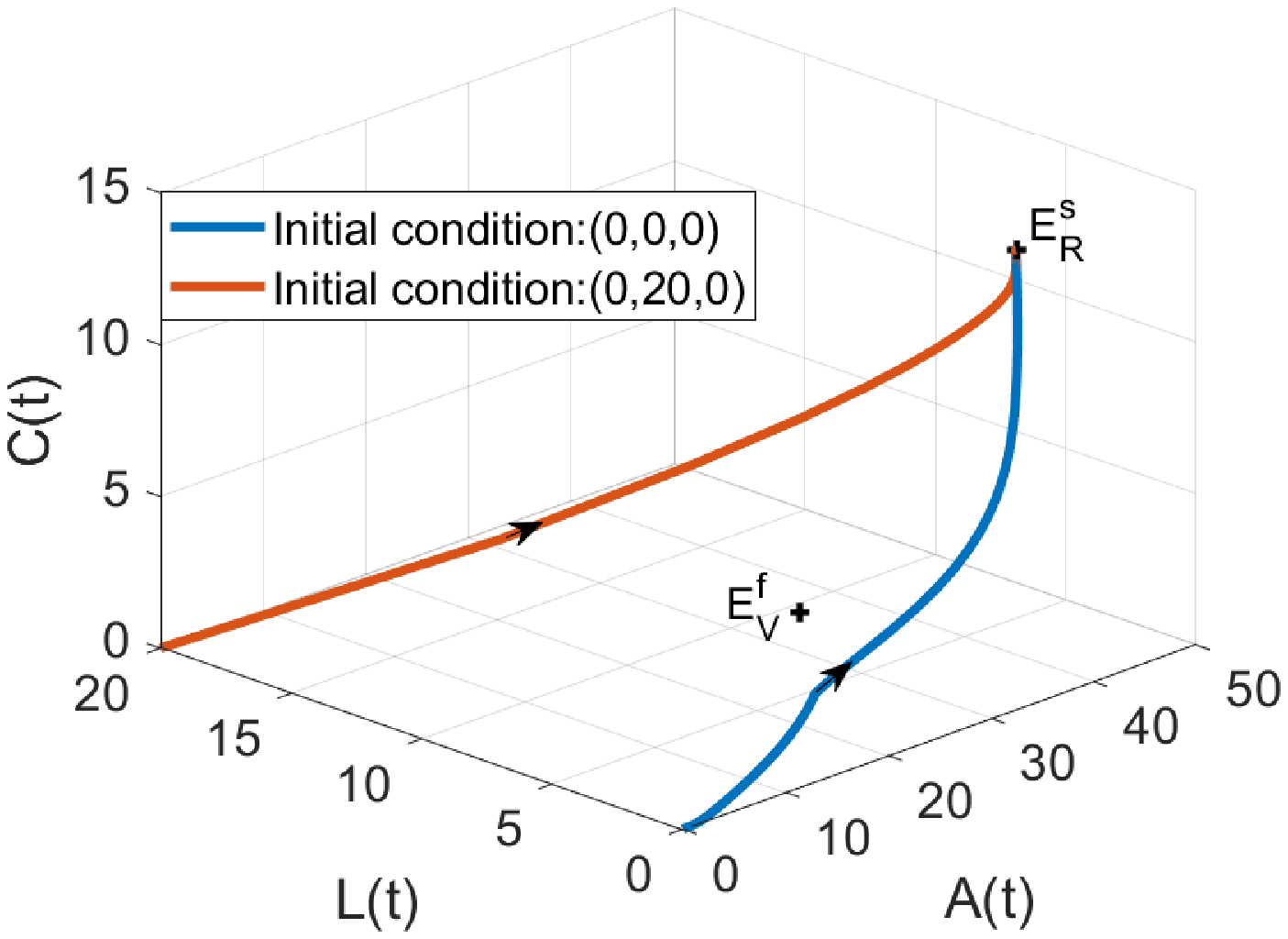}}
		\caption{Time series and phase plots of System (3.2) shows that equilibrium $E_{R}^{s}$ is a unique attractor when $N=300$ and other parameters are taken as in Figure \ref{ScenarioA} .}
		\label{ScenarioB}
	\end{figure}

	\begin{figure}[htbp]
		\centering
		\includegraphics[width=5cm]{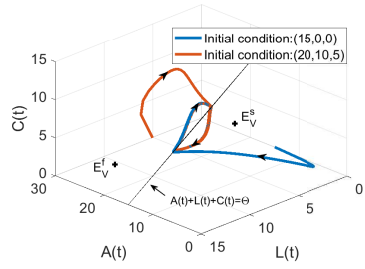}
		\caption{Phase plot of System (3.2) shows the existence of periodic solution when $E^{f}$ and $E^{s}$ are two virtual equilibria with parameters taken as in Figure 2. }
		\label{ScenarioC}
	\end{figure}
	
	\begin{figure}[htbp]
		\centering
		\subfloat[ ]{\label{TimeseriesScenarioD:subfig1}\includegraphics[width=5cm]{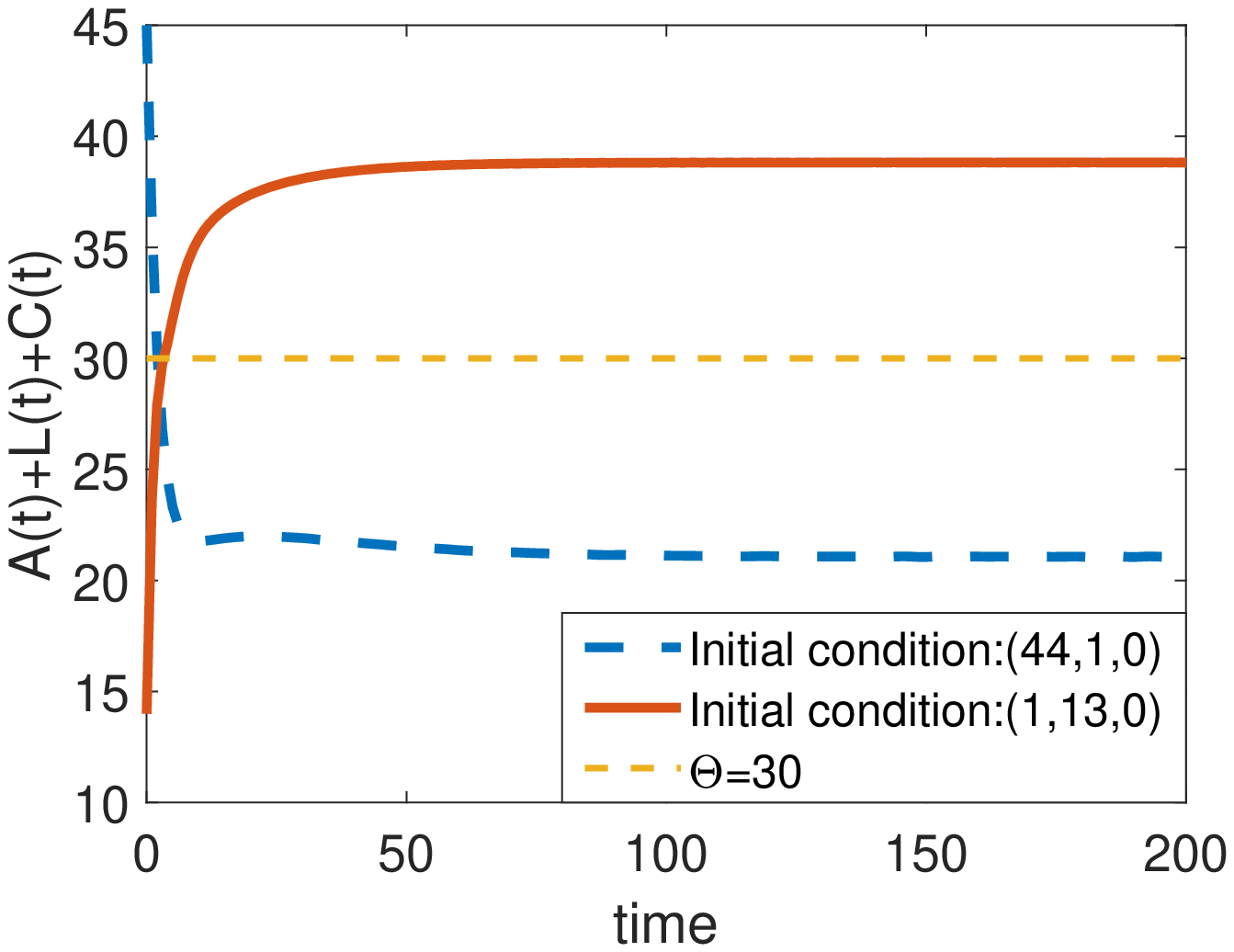}}\hspace{3mm}
		\subfloat[ ]{\label{TimeseriesScenarioD:subfig2}\includegraphics[width=5.3cm]{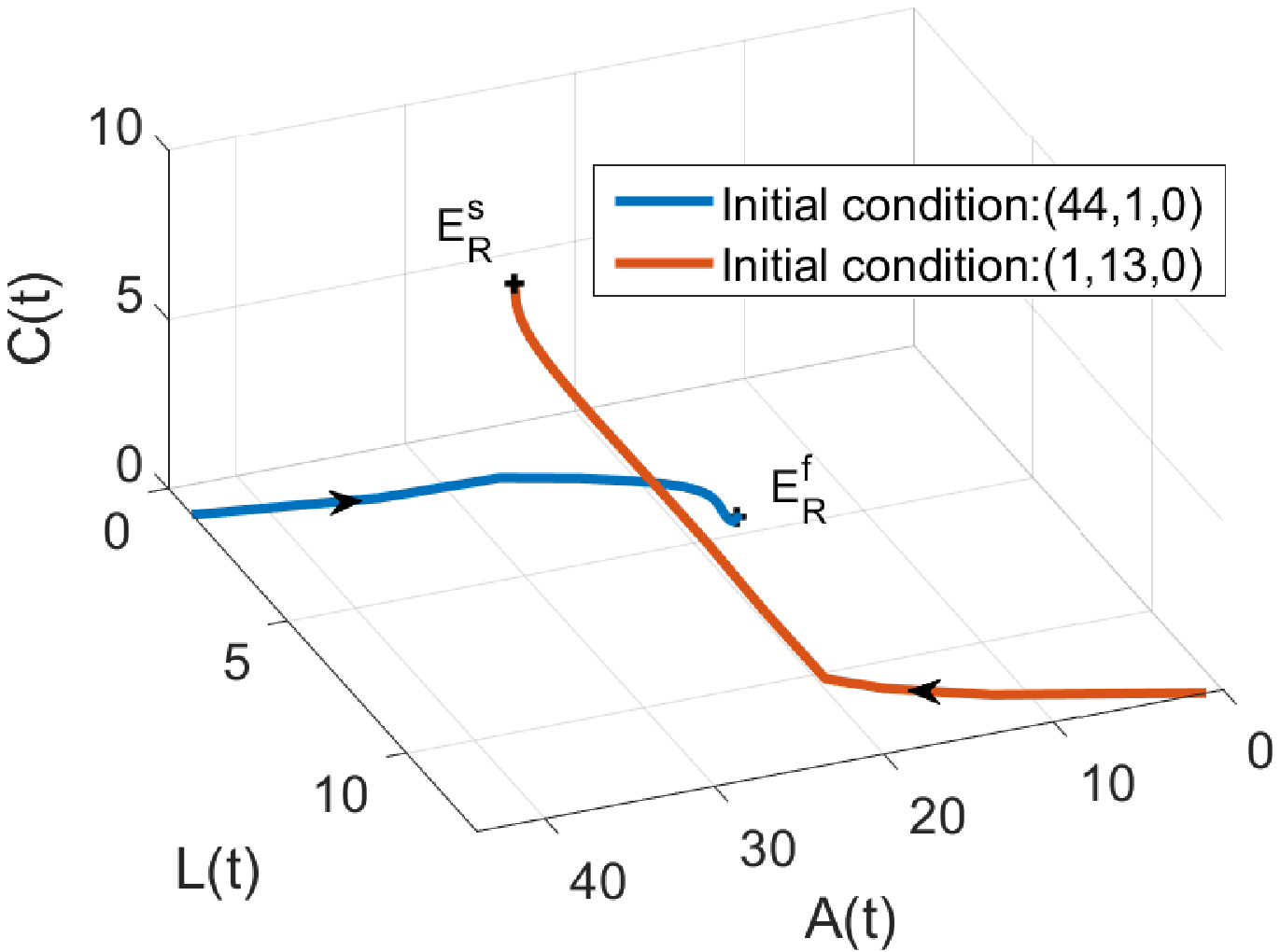}}
		\caption{Time series and phase plots of System (3.2) shows the bistablity between $E_{R}^{f}$ and $E_{R}^{s}$ when $N=200$ and other parameters are taken as in Figure \ref{ScenarioA}. }
		\label{TimeseriesScenarioD}
	\end{figure}
